\documentclass[10pt]{article}
\newcommand{\documentdate}{31 VII 2023}

\usepackage{a4wide,latexsym,amsmath,varioref,graphicx,xcolor,amssymb}
\title{Yet another fast variant of Newton's method \\ for nonconvex optimization}

\author{	
   S. Gratton%
   \thanks{Universit\'e de Toulouse, INP, IRIT, Toulouse, France. Email:
     serge.gratton@enseeiht.fr. Work partially supported by 3IA Artificial and
     Natural Intelligence Toulouse Institute (ANITI), French "Investing for the Future
     - PIA3" program under the Grant agreement ANR-19-PI3A-0004"}, 
   ~S. Jerad%
   \thanks{ANITI, Universit\'e de Toulouse, INP, IRIT, Toulouse, France. Email:
     sadok.jerad@toulouse-inp.fr}
   ~and Ph. L. Toint%
   \thanks{NAXYS, University of Namur, Namur, Belgium. Email:
     philippe.toint@unamur.be. Partly supported by ANITI.}
}

\newcommand{\beqn}[1]{\begin{equation}\label{#1}}
\newcommand{\eeqn}{\end{equation}}
\newcommand{\req}[1]{(\ref{#1})}
\newcommand{\ms}{\;\;\;\;}
\newcommand{\tim}[1]{\;\; \mbox{#1} \;\;}

\newtheorem{theorem}{Theorem}[section]
\newtheorem{lemma}[theorem]{Lemma}
\newtheorem{corollary}{Corollary}

\newcommand{\numsection}[1]{\section{#1}\setcounter{equation}{0}}
\newcommand{\appnumsection}[1]{\section*{#1}\setcounter{equation}{0}
  \renewcommand{\theequation}{A.\arabic{equation}}
  \renewcommand{\thetheorem}{A.\arabic{theorem}}
  \renewcommand{\thetable}{A.\arabic{table}}
  \renewcommand{\thefigure}{A.\arabic{figure}}
  \renewcommand{\thesection}{A} }\renewcommand{\theequation}{\arabic{section}.\arabic{equation}}

\newcounter{algo}[section]
\renewcommand{\thealgo}{\thesection.\arabic{algo}}
\newcommand{\llem}[2]{\vspace{\baselineskip} 
\noindent\framebox[\textwidth]{\parbox{0.95\textwidth}{
\begin{lemma} \label{#1} \rm #2 \end{lemma} } } \vspace{\baselineskip} }
\newcommand{\algo}[3]{\refstepcounter{algo}
\begin{center}\begin{figure}[htbp]
\framebox[\textwidth]{
\parbox{0.95\textwidth} {\vspace{\topsep}
{\bf Algorithm \thealgo : #2}\label{#1}\\
\vspace*{-\topsep} \mbox{ }\\
{#3} \vspace{\topsep} }}
\end{figure}\end{center}}
\newcommand{\bpr}{{\bf Proof.} \hspace{1.5mm}}
\newcommand{\epr}{\hfill $\Box$ \vspace*{1em}}
\newcommand{\proof}[1]{
\begin{list}{}{
\setlength{\topsep}{0.0pt}
\setlength{\partopsep}{0.0pt}
\setlength{\leftmargin}{0.025\textwidth}
\setlength{\rightmargin}{0.5\leftmargin}
\setlength{\labelwidth}{0.5\leftmargin}
\setlength{\labelsep}{0.25\leftmargin}}
\item \bpr #1 \epr \noindent
\end{list}}

\newcommand{\lthm}[2]{\vspace{\baselineskip} 
\noindent\framebox[\textwidth]{\parbox{0.95\textwidth}{
\begin{theorem} \label{#1} \rm #2 \end{theorem} } } \vspace{\baselineskip} }

\newcommand{\ii}[1]{\{ 1, \ldots, #1 \}}
\newcommand{\iiz}[1]{\{ 0, \ldots, #1 \}}
\newcommand{\iibe}[2]{\{ #1, \ldots, #2 \}}

\newcommand{\calO}{{\cal O}} 
\newcommand{\calS}{{\cal S}} 
\renewcommand{\Re}{\hbox{I\hskip -2pt R}}
\newcommand{\bigfrac}[2]{\frac{\displaystyle #1}{\displaystyle #2}}
\newcommand{\sfrac}[2]{{\scriptstyle \frac{#1}{#2}}}
\newcommand{\half}{\sfrac{1}{2}}
\newcommand{\eqdef}{\stackrel{\rm def}{=}}

\newcommand{\kap}[1]{\kappa_{\mbox{\tiny #1}}}
\newcommand{\al}[1]{{\footnotesize{\sf #1}}}
\newcommand{\tal}[1]{{\normalsize {\sf #1}}}

\newcommand{\comment}[1]{}
\newcommand{\sigkgk}{\sigma_k \| g_k\|}
\newcommand{\minlamda}{[-\lambda_{\min}(H_k)]_{\tiny +}}

\newcommand{\hg}{\widehat{g}}
\newcommand{\hH}{\widehat{H}}
\newcommand{\minlamdavp}{[-\lambda_{\min}(\hH_k)]_{\tiny +}}
\newcommand{\regstep}{\al{RegStep}}
\newcommand{\newtstep}{\al{NewtonEigenStep}}

\date{\documentdate}

\begin{document}

\maketitle

\begin{abstract}
  %\phil{
A class of second-order algorithms is proposed for minimizing smooth nonconvex
functions that alternates between regularized Newton and negative
curvature steps in an iteration-dependent subspace. In most cases, the Hessian matrix is regularized with
the square root of the current gradient and an additional term taking
moderate negative curvature into account, a negative curvature step
being taken only exceptionally. Practical variants have been detailed
where the subspaces are chosen to be the full space, or Krylov
subspaces. In the first case, the proposed
method only requires the solution of a single linear system
at nearly all iterations.  We establish that at most
$\mathcal{O}\left( |\log\epsilon|\,\epsilon^{-3/2}\right)$
evaluations of the problem's objective function and derivatives are needed
for algorithms in the new class to obtain an $\epsilon$-approximate first-order
minimizer, and at most 
$\mathcal{O}\left(|\log\epsilon|\,\epsilon^{-3}\right)$
to obtain a second-order one. Encouraging initial numerical experiments with two
full-space and two Krylov-subspaces variants are finally presented.
%} % phil
\end{abstract}


\numsection{Introduction}

It is not an understatement to say that Newton's method is a central
algorithm to solve nonlinear minimization problems, mostly because the method
exhibits a quadratic rate of convergence when close to the solution
and is affine-invariant. In the worst case, it can however be as slow
as a vanilla first-order method \cite{CartGoulToin10a}, \cite[Section~3.2]{CartGoulToin22} even when
globalized with a linesearch \cite{Nest18} or a trust region
\cite{ConnGoulToin00}. This drawback has however been circumvented by the
cubic regularization algorithm \cite{NestPoly06} and its subsequent
adaptive variants \cite{CartGoulToin11,CartGoulToin11d}, \cite[Section~3.3]{CartGoulToin22}. For nonconvex
optimization, these latter variants exhibit a worst-case
$\mathcal{O}\left(\epsilon^{-{3/2}}\right)$ complexity order to find an $\epsilon$-
first-order minimizer compared with the
$\mathcal{O}\left(\epsilon^{-2}\right)$ order of second-order trust-region
methods \cite{GratSartToin08}, \cite[Section~3.2]{CartGoulToin22}.  Adaptive cubic regularization was later
extended to handle inexact derivatives
\cite{XuRoosMaho20,YaoXuRoosMaho21,BellGuriMoriToin19,BellGuriMori21}, probabilistic models
\cite{BellGuriMori21,CartSche17}, and even schemes in which the value
of the objective function is never computed \cite{GratJeraToin22c}.

\noindent
However, as noted in \cite{Misc21}, the improvement in complexity
has been obtained by trading the simple Newton step requiring only the
solution of a single linear system for more complex or slower
procedures, such as secular iterations, possibly using Lanczos
preprocessing \cite{CartGoulToin09a,CartGoulToin11} (see also
\cite[Chapters 8 to 10]{CartGoulToin22}) or (conjugate-)gradient
descent \cite{Grie81,CarmDuch19}. In the simpler context of convex
problems, two recent papers \cite{Misc21,DoikNest23} independently
proposed another globalization technique. At an iterate $x$, the step
$s$ is computed as
\beqn{skcompute}
s = -(\nabla_x^2 f(x) + \lambda_k I_n)^{-1} \nabla_x^1 f(x)
\eeqn
where $\lambda_k \sim \sqrt{\|\nabla_x^1 f(x) \|}$. This new approach
exhibits the best complexity rate of second-order methods for convex
optimization and retains the local superlinear convergence of standard
Newton method, while showing remarkable numerical promise
\cite{Misc21}. Devising an algorithm for nonconvex functions that
can use similar ideas whenever possible appears as a natural extension.

In the nonconvex case, the Hessian may be indefinite and it is
well-known that negative curvature can be exploited to ensure progress
towards second-order points. Mixing gradient-related (possibly Newton)
and negative curvature directions has long been considered and can be
traced back to \cite{McCo77}, which initiated a line of work using
curvilinear search to find a step combining both types of
directions. The length of the step is typically tuned using an
Armijo-like condition \cite{McCo77,Gold80,MoreSore79}.  Improvements
were subsequently proposed by incorporating the curvilinear step in a
nonmonotone algorithm \cite{FerrLuciRoma96}, allowing the resolution
of large-scale problems \cite{LuciRochRoma98} or by choosing between
the two steps based on model decrease \cite{GoulLuciRomaToin00}.
Alternatively, negative curvature has also been used to regularize the
Hessian matrix, yielding the famous Goldfeld-Quandt-Trotter (GQT) method
\cite{GoldQuanTrot66}. Unfortunately, this method also involves more
complex computation to find the step and has the same global
convergence rate as first-order algorithms \cite{UedaYama14}.
The negative curvature regularization was also the subject of the more
recent paper  \cite{BiginMartinez17}, in which various Newton steps
are tried at each iteration in order to ensure the optimal ARC
\cite{CartGoulToin11,CartGoulToin11d} global rate of convergence.

One may then wonder if it is possible to devise an adaptive
second-order method using a single explictly regularized Newton step
when possible and a negative curvature direction only when necessary,
with a near-optimal complexity rate. The objective of this paper is to
show that it is indeed possible (and efficient). To this aim, we
propose a fast Newton's method that exploits negative curvature for
nonconvex optimization problems and generalizes the method proposed in
\cite{Misc21,DoikNest23} to the nonconvex case. The new algorithm
automatically adjusts the regularization parameter (without knowledge
of the Hessian's Lipschitz constant). The method either uses an
appropriately regularized Newton step taking the smallest negative
eigenvalue of the Hessian also into account or simply follows the
negative curvature otherwise. It first attempts a step along a
direction regularized by the square root of the gradient only, as
in the convex setting \cite{DoikNest23,Misc21}. In that sense, it is
inspired by the ``convex until proved guilty'' strategy advocated by
\cite{CarmDuchHindSidf17}. If this attempt fails, it obtains negative
curvature information of the Hessian, which is then used either for
regularization or to define a step along a negative-curvature
direction. In what follows, all these
operations are carried out in a specific, iteration-dependent
subspace, whose choice leads to different algoritmic variants.
We prove that these methods require at most
$\mathcal{O}\left(| \log \epsilon|\,\epsilon^{-3/2} \right)$
iterations and evaluations of the problem data to obtain an
$\epsilon$-approximate first-order critical point, which is very close
to the optimal convergence rate of second-order methods for Lipschitz
Hessian functions \cite{CartGoulToin18a}.  We also introduce an
further algorithmic variant which is guaranteed to find a
second-order critical point in at most
$\mathcal{O}\left(| \log \epsilon|\,\epsilon^{-3} \right)$ iterations.

The paper is organized as follows. Section~2 describes the general
algorithmic framework and compares it with recent work on second-order
methods. Section~3 states our assumptions and derives a bound on its
worst-case complexity for finding first-order critical points.
Section~4 presents the second-order algorithmic variant and states its
complexity, the corresponding analysis being detailed in
appendix. Section~5 then discusses some choices of the
iteration-dependent subspace, including  Krylov
spaces. Section~6 finally illustrates the numerical behavior of the
proposed methods. Conclusions are drawn in Section~7.

\noindent 
\textbf{Notation} The following notations will be used throughout the
paper. Let $q \geq 1$. The symbol $\| . \|$ denotes the Euclidean norm
for vectors in $\Re^q$ and its associated subordinate norm for matrices.
$\lambda_{\min}(M)$  and $\lambda_{\max}(M)$ denote the minimum and
maximum eigenvalues of a symmetric matrix $M$, while $I_q$ is the
identity matrix in $\Re^q$.  For $x \in \Re$, we define
$[x]_{+} = \max(x,0)$. For two vectors $x,y \in \Re^q $,
$x^\intercal y$ denotes their inner product. The $i$-th column of $I_q$ is
denoted by $e_i$.

\numsection{Adaptive Newton with Negative Curvature }\label{thealgo-s}

We consider the problem of finding approximate first-order critical
points of the smooth unconstrained nonconvex optimization problem
\beqn{minf}
\min_{x \in \Re^n} f(x)
\eeqn
and discuss our algorithm called \al{AN2C} (for
\tal{Adaptive Newton with Negative Curvature})
\vpageref{AN2C}. The algorithm, whose purpose is to compute
first-order critical points, is presented
in the framework of adaptive regularization methods 
\cite{BirgGardMartSantToin17,CartGoulToin11}
\cite[Section~3.3]{CartGoulToin22} and proceeds as follows, using
two subroutines $\regstep$ and $\newtstep$.
%Broadly speaking, the algorithm first attempts to define the step using
%a regularization depending only on a regularization parameter
%$\sigma$ and the gradient. If this fails, the 
%step is the defined using a regularization also depending on negative
%curvature. It is then accepted or rejected by comparing achieved
%to predicted reductions (as in trust-region and adaptive
%regularization methods \cite{CartGoulToin22}), and the regularization
%parameter is finally updated \cite{BirgGardMartSantToin17, CartGoulToin11}.  

\algo{AN2C}{Adaptive Newton with Negative Curvature  (\tal{AN2C})}{
	\begin{description}
		\item[Step 0: Initialization] An initial point $x_0\in \Re^n$, a regularization
		parameter $\sigma_0>0$ and a gradient accuracy threshold
		$\epsilon \in (0,1]$  are given, as well as the parameters
			\[
		\sigma_{\min}  > 0, \, \kappa_C , \kappa_\theta > 0,\,
		\kappa_a, \kappa_b \geq 1,
                \,\varsigma_1\in(0,1),\,
                \varsigma_2 \in [0,\sfrac{1}{2}),\,
                \varsigma_3 \in[0,1),\,
                \theta \in(0,1],
		\]
		\[
		0< \gamma_1 < 1 < \gamma_2 \leq \gamma_3 
		\tim{and }
		0 < \eta_1 \leq \eta_2 < 1.  
		\]
		Set $k=0$.
		\item[Step 1: Check termination] Evaluate
		${g_k} \eqdef \nabla_x^1 f(x_k)$. Terminate if
		$\|g_k \| \leq \epsilon$.
		\item [Step 2: Compute subspace derivatives]
		Choose $p \in \ii{n}$ and form $V_p \in \Re^{n \times p}$. 
		Compute $\hg_k \eqdef V_p^\intercal g_k $ and $\hH_k
                \eqdef V_p^\intercal H_k V_p$ where $H_k \eqdef \nabla_x^2 f(x_k)$.
		\item[Step 3 (Optionnal): Attempt a regularization step]
		\beqn{sdefstep}
		s_k = s_k^{def} = \regstep( \, \hg_k, \hH_k, V_p, \|g_k\|, \sigma_k, \kappa_a, \kappa_b, \varsigma_1, \varsigma_2 \, ).
		\eeqn
                If $s_k^{def}$ has been successfully defined, go to Step 5.
		\item[Step 4 : Newton Step Computation]
		\beqn{sknewton}
		s_k  = \al{NewtonEigenStep}(\, \hg_k, \hH_k, V_p, \|g_k\|, \sigma_k, \kappa_C, \kappa_b, \kappa_\theta, \varsigma_3,\theta \,).
		\eeqn
		
		\item[Step 5: Acceptance ratio computation ]
		Evaluate $f(x_k +s_k)$ and compute the acceptance ratio
		\beqn{rhokdef}
		\rho_k = \frac{f(x_k) - f(x_k +s_k)}{-(g_k^\intercal s_k + \frac{1}{2} s_k^\intercal H_k s_k)}.
		\eeqn
		If $\rho_k \geq \eta_1$, set $x_{k+1} = x_k + s_k$ else $x_{k+1} = x_k$.
		
		\item[Step 6: Regularization parameter update] Set
		\begin{equation}\label{sigmakupdate}
		\sigma_{k+1} \in \left\{
		\begin{aligned}
		& [ \max\left(\sigma_{\min}, \gamma_1\sigma_k \right),\sigma_k] &&\text{ if } \rho_k \geq \eta_2, \\
		& [\sigma_k ,\gamma_2\sigma_k ]            &&\text{ if } \rho_k \in [\eta_1, \eta_2),\\
		& [\gamma_2 \sigma_k,\gamma_3\sigma_k]     &&\text{ if } \rho_k < \eta_1.
		\end{aligned}
		\right.
		\end{equation}
		Increment $k$ by one and go to Step~1.
	\end{description}
}
                  
\algo{regstep}{\regstep$(\,\hg_k, \hH_k, V_p, \|g_k\|, \sigma_k, \kappa_a, \kappa_b, \varsigma_1, \varsigma_2\,)$}{
	Attempt to solve the linear system 
	\beqn{convexstep}
	(\hH_k +  \sqrt{ \kappa_a \sigma_k \|g_k \|} I_p) y_k^{def} = - \hg_k.
	\eeqn
        If a solution $y_k^{def}$ of this system can be obtained such
        that
	\begin{align}
	(y_k^{def})^\intercal( \hH_k &+  \sqrt{\kappa_a \sigma_k \|g_k \|} I_p) y_k^{def} > 0, \label{possystem} \\
	\|y_k^{def} \| &\leq \frac{(1+\kappa_\theta)}{\varsigma_1} \sqrt{\frac{  \|g_k \|}{ \kappa_a \sigma_k}},
	\label{skconvexbound} \\
	\| H_k V_p y_k^{def} + g_k \| &\leq \kappa_b \| \hH_k y_k^{def} + \hg_k \|, \label{Vpspacequalitydef} \\
	\|r_k^{def} \|
	&\leq  \min \left( \varsigma_2 \sqrt{ \kappa_a \sigma_k \|g_k \|}\|y_k^{def} \|, \kappa_\theta \| \hg_k \| \right) \label{residcondconv} 
	\end{align}
        where $r_k^{def} = (\hH_k +  \sqrt{\kappa_a \sigma_k \|g_k \|} I_p) y_k^{def} + \hg_k$,
	then return $ s_k^{def} \eqdef V_p y_k^{def}$.
}

\algo{NewtonEigenStep}{NewtonEigenStep$(\,\hg_k, \hH_k, V_p, \|g_k\|,
  \sigma_k  , \kappa_C, \kappa_b, \kappa_\theta, \varsigma_3, \theta \,)$}{
	\begin{description}
	\item[Step 1: Test negative curvature]
	If $\lambda_{\min}(\hH_k) \leq - \kappa_C \sqrt{\sigma_k \|g_k \|}$, go to Step 4.
	
	\item[Step 2: Newton Step] Solve
	\begin{quote}
		\vspace*{-1em}
		\beqn{steplargegrad}
		\left( \hH_k + (\sqrt{\sigma_k \| g_k\|} + [-\lambda_{\min}(\hH_k)]_+) I_p \right) y_k^{neig}   = - \hg_k
		\eeqn
		to ensure the residual condition
		\beqn{residcond}
		\begin{array}{lcl}
	        \|r^{neig}_k \|
		& \eqdef&\left\|\left(\hH_k+(\sqrt{\sigma_k\| g_k\|}+[-\lambda_{\min}(\hH_k)]_+)I_p\right) y_k^{neig}+\hg_k\right\|\\
		& \leq &  \min \left( \varsigma_3 \sqrt{\sigma_k \|g_k\|}\|  y_k^{neig} \|, \kappa_\theta \| \hg_k \| \right).
		\end{array}
		\eeqn
	\end{quote}
	\item[Step 3: Check global quality of the solution] If 
	\beqn{Vpspacequality}
	\| H_k  V_p  y_k^{neig} + g_k \|
        \leq \kappa_b \| \hH_k y_k^{neig} + \hg_k \| \tim{ then set } s_k = s_k^{neig} \eqdef V_p  y_k^{neig}.
	\eeqn
	Else, go back \al{AN2C}[Step 2].
	
	\item[Step 4: Eigenvector direction] 
		Compute $u_k$ such that 
		\beqn{negcurvvector}
		\hg_k^\intercal u_k  \leq 0, \, \,  \|u_k \| = 1
		\tim{and}
                u_k^\intercal \hH_k u_k \leq \theta \lambda_{\min}(\hH_k)
		\eeqn
		and set 
		\beqn{negcurvstep}
		s_k = s_k^{curv} \eqdef  \frac{\theta \kappa_C \sqrt{\sigma_k \|g_k \|}}{\sigma_k} V_p u_k.
		\eeqn
	\end{description}

}

The selection of the iteration-dependent subspace defined as the
range of $V_p$ in Step~2 is of course crucial for the algorithm.
At this stage of the algorithm description, cycling may possibly occur
between Step~2 and \req{Vpspacequality} in Step~4, should the choice of the subspace be
consistently inadequate. We will however  discuss some practical
choices in Section~\ref{numericstep}, for which this situation cannot
happen. For our subsequent analysis, we therefore assume the following.

\noindent
\textbf{AS.0} For each iteration $k$, condition \req{Vpspacequality} is satisfied
after finitely many choices of $V_p$.  Moreover, there exists a constant $V_{\max}\ge 1$ such that
\beqn{Vpbound}
\|V_p\| \leq V_{\max} \tim{for all } p \in \ii{n}.
\eeqn 

\noindent
After selecting the subspace\footnote{Since we do not specify at this
point how to make this selection, \al{AN2C} may be viewed as a \emph{class of algorithms}
depending on the choice of $V_p$.} and projecting the current gradient and
Hessian, we first attempt a step that avoids computing negative
curvature information.  Indeed, the $s_k^{def}$ notation, where $def$
stands for ``definite'', in \req{convexstep} makes the connection with
the two conditions \eqref{convexstep} and \eqref{possystem}. The
condition \eqref{possystem} is significantly less restrictive than
checking the positive-definess of the regularized matrix in
\eqref{convexstep}. This is at variance with the work of
\cite{BirginMartinez19} where a factorization is required at each
step, and coherent with the 'capped-CG' subroutine proposed at
\cite{YaoXuRoosMahoWri22}. Should the problem be (locally) convex,
\req{possystem} would automatically hold (see
\cite{Misc21,DoikNest23}).  The test \eqref{skconvexbound} is required
as to avoid steps whose magnitude is too large compared to the
gradient (the motivation for its particular form of the test will
become clear in Section~\ref{complexity-s}).

When computing a vector 
$s_k^{def}$ satisfying \req{possystem}  to
\req{residcondconv} is
not possible, we (approximately) solve a linear system in $\Re^{p}$
\eqref{steplargegrad} whose definition involves $\minlamdavp$.
Even if an exact solution can be obtained at a marginal cost for small
$p$, we still allow an approximate solution satisfying
\eqref{residcond}. We note that $\minlamdavp$ could have been
replaced in  \eqref{steplargegrad} by $\kappa_C \sqrt{\sigkgk}$  and the
remainder of the complexity analysis would remain valid.  
An interesting connection can also be established between the
regularization in \req{steplargegrad} for $\hH_k= H_k$ and the GQT method
\cite{GoldQuanTrot66}, as the regularization parameter
$(\sqrt{\sigma_k \| g_k\|} + [-\lambda_{\min}(H_k)]_+  )$ is very
similar in spirit to that used in this method.
In the closely related algorithm  of \cite{BiginMartinez17}, a term
$\mu$ is added to $ [-\lambda_{\min}(H_k)]_+ $ and multiple $\mu$'s
are tested as to ensure 'cubic' descent. In our case, $\sqrt{\sigkgk}$
directly yields a regularization of the desired order. 
Also observe that, in most cases, the ``approximate minimum curvature
direction'' $u_k$ is already available when
computing $\lambda_{\min}(\hH_k)$. It can be also retrieved via a
Lanczos procedure as proposed in \cite[Lemma 9]{RoyeWrig18}. 

We now provide some comments that apply to the definition of
both $s_k^{def}$ and $s_k^{neig}$. Specifically, focusing on the
latter, condition \eqref{Vpspacequality} serves to ensure the
appropriateness of the subspace spanned by $V_p$. This condition
guarantees that the projected residual \req{residcond} is sufficiently
small compared to both the projected and unprojected gradients. In a
more standard setting, where $V_p = I_n$, this condition simplifies to 
\beqn{residallcase}
\|(H_k + (  \sqrt{\sigkgk} + \minlamda   ) I_n ) s_k + g_k \|
\leq \min (\varsigma_3 \sqrt{\sigkgk} \|s_k \|, \kappa_\theta \|g_k \|),
\eeqn
where the  $\kappa_\theta \|g_k \|$ term is standard when devising
truncated CG algorithms. The other term ensures the typical condition
required for the approximate minimization of the cubic model $m_k$,
namely that 
\beqn{mkcubiccond}
\| \nabla_s^1 m_k(s_k) \| \leq \mathcal{O}(\|s_k \|^2).
\eeqn
This condition is typically used to derive the optimal complexity rate
$\mathcal{O} \left( \epsilon^{-3/2}\right)$, see
\cite{CartGoulToin11,BirgGardMartSantToin17} and the references
therein.

Finally we note that the condition of Step~1 in the
\al{NewtonEigenStep} algorithm, which forces the negative curvature
step \req{negcurvstep},  can be interpreted as the comparison of the
minimal curvature of the quadratic ($\lambda_{\min}(V_p^\intercal H_k
V_p)$) with the quantity $\sigma_k \sqrt{\|g_k\|/\sigma_k}$, which
itself can be viewed as the curvature of the regularization term
$\sfrac{1}{6}\sigma_k\|s\|^3$ for some $s$ whose length
$\sqrt{\|g_k\|/\sigma_k}$ is of the order of a standard regularized
step (see \cite[Lemma~3.3]{CartGoulToin22}, for instance). The test
thus ensures a ``regularization-like'' step when the quadratic's
negative curvature is strong enough to dominate that of the
regularization for too small steps (see \req{lowboundcurv} below).
    
Once the step has been computed, the mechanisms for
accepting/rejecting the new iterate (Step~5) and updating the regularization
parameter (Step~6) are typical of adaptive regularization algorithms
(see \cite{BirgGardMartSantToin17, CartGoulToin11}
or \cite[Section~3.3.1]{CartGoulToin22}, for instance).

Before delving into the complexity analysis of \al{AN2C}, we further
explore its fundamental properties and discuss its relationships with
closely related nonconvex optimization algorithms. The method
presented in \cite{CurtRobi18} differs from \al{AN2C} in that it
employs a gradient step followed by a negative curvature step. On the
other hand, \cite{LiuLiWangYiYang18} adopts a condition-based approach
to choose between gradient descent and negative curvature directions,
relying on known smoothness parameters, while our methods remain fully
adaptive. Another related approach is presented in \cite{RoyeWrig18},
which, unlike \al{AN2C}, examines various conditions to select a
specific direction (gradient, Newton, negative curvature) and performs
a linesearch. Furthermore, \cite{CurtRobiRoyeWrig21} proposes a
trust-region algorithm (in contrast to adaptive regularization) that
tackles the trust-region subproblem using a combination of conjugate
gradients and negative curvature. Notably, their condition on the
residuals of this subproblem
\cite[Inequality~(3.2)]{CurtRobiRoyeWrig21} can be related to
\eqref{residallcase}. 

Following well-established practice, we now define
\[
\calS \eqdef \{ k \geq 0 \mid x_{k+1} = x_k+s_k \} = \{ k \geq 0 \mid \rho_k \geq \eta_1 \},
\]
the set of indexes of ``successful iterations'', and
\[
\calS_k \eqdef \calS \cap \iiz{k},
\]
the set of indexes of successful iterations up to iteration $k$. We
further partition $\calS_k$ in three subsets depending on the nature of
the step taken, so that
\[
\calS_k^{neig} \eqdef \calS_k \cap \{s_k = s_k^{neig} \},
\quad \calS_k^{curv} \eqdef \calS_k \cap \{s_k = s_k^{curv} \},
\quad \calS_k^{def} \eqdef \calS_k \cap \{s_k = s_k^{def} \}.
\]
We also recall a well-known result bounding the total number of
iterations of adaptive regularization methods in
terms of the number of successful ones.

\llem{SvsU}{\cite[Lemma~2.4]{BirgGardMartSantToin17}, \cite[Lemma~2.4.1]{CartGoulToin22}
	Suppose that the \al{AN2C} algorithm is used and that $\sigma_k \leq
	\sigma_{\max}$ for some $\sigma_{\max} >0$. Then
	\beqn{unsucc-neg}
	k \leq |\calS_k| \left(1+\frac{|\log\gamma_1|}{\log\gamma_2}\right)+
	\frac{1}{\log\gamma_2}\log\left(\frac{\sigma_{\max}}{\sigma_0}\right).
	\eeqn
}
\proof{Since Steps~5 and 6 of the \al{AN2C} algorithm are identical to
  Steps~3 and 4 of the algorithm presented on page~43 of
  \cite{CartGoulToin22} and since Lemma~2.4.1 in this reference only
  depends on these steps, the conclusion of the lemma remain valid for
  \al{AN2C} algorithm. } 

This result implies that the overall complexity of the algorithm can be
estimated once bounds on $\sigma_k$ and $|\calS_k|$ are known, as we will show in the
next section.

We now state a simple relation between
 $\|s_k^{neig} \|$, $\sigma_{k}$ and $\| g_k\|$ inspired by
 \cite{Misc21}.

\llem{skstepbound}{
For all iterations $k$ where $s_k^{neig}$ is computed, we have that
\beqn{skneqyp}
\|s_k^{neig} \| \leq V_{\max} \|  y_{k}^{neig} \|
\eeqn

\beqn{gkexpressionnewton}
\hg_k = -\left(\hH_k+ (\sqrt{\sigma_k \| g_k\|} + [-\lambda_{\min}(\hH_k)]_+  ) I_p      \right)  y_{k}^{neig} + r_k^{neig}
\eeqn
and
\beqn{sknnormbound}
\|s_k^{neig} \| \leq (1+\kappa_\theta) V_{\max}^\sfrac{3}{2} \sqrt{\frac{\|g_k \|}{\sigma_k}}.
\eeqn
Similarly, when $s_k^{def}$ is computed,
\beqn{gkconvexpression}
\hg_k = -\left(\hH_k + \sqrt{ \kappa_a \sigkgk} I_p\right)  y_{k}^{def} + r_k^{def},
\eeqn  
and
\beqn{skdefbound}
\|s_k^{def} \| \leq V_{\max} \| y_{k}^{def} \| .
\eeqn
At last, when $s_k^{curv}$ is computed,
\beqn{skcurvbound}
\| s_k^{curv}\| \leq V_{\max} \frac{\theta \kappa_C \sqrt{\sigkgk}}{\sigma_k}. 
\eeqn
}
\proof{
First note from the second part of \eqref{Vpspacequality} and \eqref{Vpbound}, 
\[
\|s_k^{neig} \| = \| V_p  y_{k}^{neig} \| \leq V_{\max} \| y_{k}^{neig} \|
\]	
yielding \eqref{skneqyp}. A similar proof can be followed to derive \eqref{skdefbound}.
	
Equation \eqref{gkexpressionnewton} results from \eqref{steplargegrad}, \eqref{Vpspacequality}
and the definition of the residual \eqref{residcond}.
Let us rewrite now \eqref{gkexpressionnewton} in function of $\hg_k$ and $\hH_k$,
\[
\hg_k = -\left(\hH_k + (\sqrt{\sigma_k \| g_k\|} + [-\lambda_{\min}(\hH_k)]_+  ) I_p \right) y_{k}^{neig} + r_k^{neig}.
\]
 
From \eqref{gkexpressionnewton}, the facts that
$\hH_k + (\sqrt{\sigkgk} + [-\lambda_{\min}(\hH_k)]_+) I_p$
is a positive definite matrix with 
\[
\lambda_{\min}(\hH_k + (\sqrt{\sigkgk} + [-\lambda_{\min}(\hH_k)]_+) I_p) \geq \sqrt{\sigkgk}
\]
and that $\|r_k^{neig} \| \leq \kappa_\theta \|\hg_k \|$ because of
\eqref{residcond}. We thus obtain that
\beqn{ypbound}
\| y_{k}^{neig} \|
\leq (1+\kappa_\theta) \sqrt{\frac{\| \hg_k \|}{\sigma_k}}
\leq (1+\kappa_\theta) \sqrt{\frac{ V_{\max}\|g_k \|}{\sigma_k}},
\eeqn
where the last inequality follows from \eqref{Vpbound}. This last
inequality and \eqref{skneqyp} give \eqref{sknnormbound}.

If $k \in \calS_k^{def}$, \eqref{gkconvexpression} is obtained from
\eqref{convexstep} and the definition of the residual $r_k^{def}$.

Else if $k \in \calS_{k}^{curv}$, \eqref{Vpbound},  the fact that
$\| u_k \| = 1$ and \eqref{negcurvstep} give \eqref{skcurvbound}.
}

\noindent
The next lemma gives a lower bound on the decrease of the local
quadratic approximation. In standard adaptive regularization algorithms, this
decrease automatically results from the minimization of the model
(See \cite{BirgGardMartSantToin17} for instance). In our case, we need to use the
properties of $s_k^{def}$, $s_k^{curv}$ and $s_k^{neig}$ to obtain
the desired result.  

\llem{lowratio}{ Let $k$ be a successful an iteration of \al{AN2C}.
If $k \in \calS_k^{def}$, we have that
\beqn{lowboundconv}
-\left(g_k^\intercal s_k + \frac{1}{2} s_k^\intercal H_k s_k \right)
\geq \frac{1-2\varsigma_2 }{2} \sqrt{ \kappa_a \sigkgk} \| y_{k}^{def} \|^2 \geq \frac{1-2\varsigma_2 }{2} \sqrt{ \kappa_a \sigkgk} \frac{\|s_k \|^2}{V_{\max}^2}.
\eeqn	
If $k \in \calS_k^{neig}$, then
\beqn{lowboundnewton}
-\left(g_k^\intercal s_k + \frac{1}{2} s_k^\intercal H_k s_k \right)\geq (1-\varsigma_3) \sqrt{  \sigkgk} \| y_{k}^{neig} \|^2 \geq
(1-\varsigma_3) \sqrt{\sigkgk}\frac{ \|s_k \|^2}{V_{\max}^2}.
\eeqn
Else, if $k \in \calS_k^{curv}$,
\beqn{lowboundcurv}
-\left(g_k^\intercal s_k + \frac{1}{2} s_k^\intercal H_k s_k \right) \geq \frac{1}{2} \theta^3 \kappa_C^3 \frac{\|g_k\|^\sfrac{3}{2}}{\sqrt{\sigma_k}} \geq \frac{1}{2} \sigma_{k} \frac{\|s_k \|^3}{V_{\max}^3}.
\eeqn
}
\proof{
Suppose first that $k \in \calS_k^{def}$. We then obtain from
\eqref{sdefstep}, \eqref{residcondconv} and \eqref{skconvexbound} that
\begin{align*}
g_k^\intercal s_k^{def} + \frac{1}{2} (s_k^{def})^\intercal H_k s_k^{def} &= ( V_p^\intercal g_k)^\intercal  y_{k}^{def} + \frac{1}{2} ( y_{k}^{def})^\intercal V_p^\intercal H_k V_p  y_{k}^{def} \\  
&=  (r_k^{def})^\intercal  y_{k}^{def} - ( y_{k}^{def})^\intercal (\hH_k
+ \sqrt{ \kappa_a \sigkgk} I_p)  y_{k}^{def} \\
&\hspace*{5mm} + \frac{1}{2} ( y_{k}^{def})^\intercal \hH_k  y_{k}^{def} \\
&= -\sqrt{\kappa_a\sigkgk}\| y_{k}^{def} \|^2+(r_k^{def})^\intercal y_{k}^{def}
      -\frac{1}{2}( y_{k}^{def})^\intercal \hH_k   y_{k}^{def} \\
&\leq - \sqrt{ \kappa_a \sigkgk} \| y_{k}^{def} \|^2
      + \varsigma_2 \sqrt{ \kappa_a \sigkgk} \| y_{k}^{def} \|^2 \\
&\hspace*{5mm} + \frac{1}{2} \sqrt{ \kappa_a \sigkgk} \| y_{k}^{def} \|^2. 
\end{align*}
Hence \eqref{lowboundconv}  follows  froms \eqref{skdefbound}.  
	
Suppose now that $k \in \calS_k^{neig}$. By using
\eqref{gkexpressionnewton} and  the fact that $ \hH_k + \minlamdavp I_p
\succeq 0$,
\begin{align*}
  &g_k^\intercal s_k^{neig} + \frac{1}{2} (s_k^{neig})^\intercal H_k s_k^{neig} =  (V_p^\intercal g_k)^\intercal  y_{k}^{neig} + \frac{1}{2} ( y_{k}^{neig})^\intercal V_p^\intercal H_k V_p  y_k^{neig} \\
&= (r_k^{neig})^\intercal  y_{k}^{neig} - ( y_{k}^{neig})^\intercal ( \hH_k + \minlamdavp I_p)  y_{k}^{neig} \\
&\hspace*{5mm} + \frac{1}{2} ( y_{k}^{neig})^\intercal (\hH_k +  \minlamdavp I_p)  y_{k}^{neig} \\
&\hspace*{5mm}- \frac{1}{2} \minlamdavp  \| y_{k}^{neig} \|^2 -\sqrt{\sigkgk} \| y_{k}^{neig} \|^2 \\
&= (r_k^{neig})^\intercal   y_{k}^{neig}  - \frac{1}{2} ( y_{k}^{neig})^\intercal (\hH_k +  \minlamdavp I_p)  y_{k}^{neig}  \\
&\hspace*{5mm}- \frac{1}{2} \minlamdavp  \| y_{k}^{neig} \|^2 -\sqrt{\sigkgk} \| y_{k}^{neig} \|^2 \\
&\leq \varsigma_3 \sqrt{\sigkgk} \|  y_{k}^{neig}\|^2
      -\frac{1}{2} \minlamdavp  \| y_{k}^{neig} \|^2 -\sqrt{\sigkgk} \| y_{k}^{neig}\|^2,
\end{align*}
where we have used \eqref{residcond} to obtain the last inequality. 
Rearranging, ignoring the $\frac{1}{2} \minlamdavp  \| y_{k}^{neig}\|^2$ term and  using \eqref{skneqyp}  yield \eqref{lowboundnewton}.

Suppose finally that $k \in \calS_k^{curv}$. As \eqref{negcurvvector}
and \eqref{negcurvstep} hold and that Step~4 of \al{NewtonEigenStep} is taken when $\lambda_{\min}(\hH_k) \leq - \kappa_C \sqrt{\sigkgk}$, we deduce that
\begin{align}
g_k^\intercal s_k^{curv}+\frac{1}{2}(s_k^{curv})^\intercal H_k s_k^{curv}
&= g_k^\intercal V_p u_k +\frac{1}{2}(s_k^{curv})^\intercal H_k s_k^{curv} \nonumber \\
&\leq \frac{1}{2} \frac{\theta^2 \kappa_C^2 \|g_k\|}{\sigma_k} u_k^\intercal \hH_k u_k
\leq \frac{1}{2} \frac{\theta^3 \kappa_C^2 \|g_k\|}{\sigma_k}  \lambda_{\min}(\hH_k) \nonumber   \\
&\leq - \frac{1}{2} \theta^3 \kappa_C^3 \frac{\|g_k\|^\sfrac{3}{2}}{\sqrt{\sigma_k}}, \label{curvquad}
\end{align}

yielding the first inequality in \eqref{lowboundcurv}. For the second inequality, remark that from \eqref{skcurvbound}, we derive that
\[
\theta^3 \kappa_C^3 \frac{\|g_k\|^\sfrac{3}{2}}{\sqrt{\sigma_k}} = \sigma_k \theta^3 \kappa_C^3 \frac{\|g_k\|^\sfrac{3}{2}}{\sigma_{k}^\sfrac{3}{2}} \geq \sigma_k \frac{\|s_k\|^3}{V_{\max}^3},  
\]
injecting the last bound in \eqref{curvquad} gives the second inequality in \eqref{lowboundcurv}.
}
\numsection{Complexity analysis for the \tal{AN2C} algorithm}\label{complexity-s}

We now turn to analyzing the worst-case complexity of the \al{AN2C}
algorithm. Our analysis is
conducted under AS.0 and the following assumptions. 

\noindent
\textbf{AS.1} The function $f$ is two times continuously differentiable in $\Re^n$. 

\noindent
{\bf AS.2} There exists a constant $f_{\rm low}$ such that
$f(x) \geq f_{\rm low}$ for all $x\in \Re^n$.

\noindent
\textbf{AS.3} The  Hessian of $f$ is globally Lipschitz
continuous, that is, there exists a non-negative constant $L_H$ such that
\beqn{LipHessian}
\|\nabla_x^2 f(x) - \nabla_x^2 f(y)\| \leq L_H \|x-y\| \, \text{ for all } x,y \in \Re^n.
\eeqn

\noindent
\textbf{AS.4}
There exists a constant $\kappa_B > 0$ such that
\[
\max(0,-\lambda_{\min}(\nabla_x^2 f(x))) \leq \kappa_B
\tim{for all} x\in \{ y \in \Re^n \mid f(y)\leq f(x_0)\}.
\]
AS.1-AS.3 are standard assumptions when analyzing algorithms that
utilize second-order information
\cite{CartGoulToin11d,BirgGardMartSantToin17}. AS.4 is weaker than
assuming bounded Hessians, a condition often used
when theoretically analyzing second-order methods that combines
negative curvature and gradient based directions
\cite{RoyeWrig18,CurtRobi18,LiuLiWangYiYang18}.
The left-hand side of the inequality is sometimes called the ``convex
deviation'' or ``modulus of nonconvexity'' \cite{KongLewi22}.
As it turns out, AS.4 is only needed for $x$ being any
iterate $x_k$ produced by the algorithm and these iterates all
belong to the level associated with the starting point $x_0$ because
the acceptance condition in Step~5 ensures that the sequence $\{f(x_k)\}$ is
non-increasing. If this level set is bounded or if the sequence
$\{x_k\}$ remains bounded for any other reason, we immediately obtain that
\beqn{boundnegcurv}
\max(0,-\lambda_{\min}(H_k))\leq \kappa_B \tim{ for all } k\geq 0
\eeqn
for some $\kappa_B \geq 0$, and both AS.3 and AS.4 automatically hold.

Having established a lower bound on the decrease ratio in
Lemma~\ref{lowratio}, we next proceed to derive an upper bound on the
regularization parameter. This is a crucial step when analyzing
adaptive regularization methods.

\llem{boundsigmak}{
Suppose that AS.1 and AS.3 hold. Then, for all $k \geq0$,
\beqn{sigmamax}
\sigma_{k} \leq \sigma_{\max} \eqdef \gamma_3 \max \left( \sigma_0, \, \varsigma_{\max} \frac{L_H}{6(1-\eta_2)} \right),
\eeqn
where 
\beqn{varsigmamaxdef}
\varsigma_{\max} \eqdef \max \left(\frac{(1+\kappa_\theta) V_{\max}^\sfrac{7}{2}}{(1-\varsigma_3)},
\, \frac{2(1+\kappa_\theta) V_{\max}^3}{\kappa_a \varsigma_1 (1-2\varsigma_2 )}, \, 2\, V_{\max}^3\right).
\eeqn
}
\proof{
Let us compute the ratio $\rho_k$ for $k\in \calS_k^{neig}$. By using
AS.3 and the standard error bound for Lispschitz approximation of the
function (see \cite[Lemma~2.1]{CartGoulToin20b}), that $\varsigma_3 < 1$,
\eqref{lowboundnewton} and \eqref{sknnormbound}, we obtain that
\begin{align}\label{rhoskn}
1-\rho_k
& = \frac{f(x_k + s_k) - f(x_k) -g_k^\intercal s_k - \frac{1}{2} s_k^\intercal H_k s_k }
       {-(g_k^\intercal s_k +\frac{1}{2} s_k^\intercal H_k s_k)} \nonumber\\
&\leq \frac{L_H V_{\max}^2 \|s_k^{neig} \|^3}{6 (1-\varsigma_3)\sqrt{\sigkgk} \|s_k^{neig} \|^2} \nonumber \\ 
       &\leq \frac{L_H V_{\max}^2 \|s_k^{neig} \|}{6 (1-\varsigma_3) \sqrt{\sigkgk} } \nonumber\\
&\leq \frac{L_H(1+\kappa_\theta) V_{\max}^\sfrac{7}{2}}{6 (1-\varsigma_3) \sigma_{k}} .
\end{align}
Hence, if
$\sigma_{k} \geq \frac{L_H(1+\kappa_\theta) V_{\max}^\sfrac{7}{2}}{6 (1-\varsigma_3)(1-\eta_2)} $,
then $\rho_k \geq \eta_2$, which implies that iteration $k$ is
successful and $\sigma_{k+1}\leq \sigma_{k}$ because of \eqref{sigmakupdate}. The mechanism of
\req{sigmakupdate} in the algorithm then ensures that 
\beqn{firstmax}
\sigma_{k} \leq \gamma_3 \max \left( \sigma_0, \frac{L_H (1+\kappa_\theta) V_{\max}^\sfrac{7}{2}}{6 (1-\varsigma_3)(1-\eta_2)} \right).
\eeqn

Similarly, if $k \in \calS_k^{def}$, we use AS.3, the Lipschitz
approximation error bound, the fact that $\varsigma_2 < \frac{1}{2}$, \eqref{lowboundconv}, \eqref{skconvexbound} and \eqref{skdefbound} to
deduce that
\[
1-\rho_k
\leq \frac{L_H \|s_k^{def} \| V_{\max}^2}{3 (1-2\varsigma_2 ) \sqrt{ \kappa_a \sigkgk} }
\leq  \frac{ L_H (1+\kappa_\theta) V_{\max}^3}{3 \kappa_a \varsigma_1 (1-2\varsigma_2 )  \sigma_{k}}.
\]
Using the same argument as above, we now obtain that
\beqn{secmax}
\sigma_{k}
\leq \gamma_3 \max\left(\sigma_0,\frac{L_H(1+\kappa_\theta) V_{\max}^3}
           {3\kappa_a\varsigma_1(1-2\varsigma_2)(1-\eta_2)}\right).
\eeqn
Consider finally the case where $k \in \calS_k^{curv}$. Again using
AS.3, the Lipschitz approximation error bound and \eqref{lowboundcurv}
lower-bound, we derive that
 \[
 1-\rho_k
 = \frac{f(x_k+s_k)-f(x_k)-g_k^\intercal s_k-\frac{1}{2}s_k^\intercal H_ks_k}
        {-g_k^\intercal s_k - \frac{1}{2} s_k^\intercal H_k s_k}
\leq \frac{L_H\|s_k^{curv}\|^3 V_{\max}^3}{6 \frac{1}{2}\sigma_k\|s_k^{curv}\|^3}
= \frac{L_H V_{\max}^3}{3 \sigma_{k}},
 \]
so that
\beqn{thirdmax}
\sigma_{k} \leq \gamma_3 \max \left( \sigma_0, \frac{L_H V_{\max}^3}{3(1-\eta_2)} \right).
\eeqn
Combining \eqref{firstmax}, \eqref{secmax} and \eqref{thirdmax} gives
\eqref{sigmamax} with $\varsigma_{\max}$ defined by \eqref{varsigmamaxdef}.
}
We now prove a lower bound on the decrease at a
successful iteration $k$ using negative curvature. We will also bound
the change in the norm $\|g_{k+1} \|$ in term of $\|g_k \|$, which
will be useful later to bound the cardinal of a subset of
$\calS_k^{neig} \cup \calS_k^{curv}$. 

\llem{negcurb}{Suppose that AS.1, AS.3 and AS.4 hold and that $k \in \calS_k^{curv}$ before termination. Then
\beqn{threehalfneg}
f(x_k) - f(x_{k+1}) \geq  \frac{\eta_1 \theta^3 \kappa_C^3}{2\sqrt{\sigma_{\max}}}\, \epsilon^\sfrac{3}{2},
\eeqn
and 
\beqn{succgradnextineq}
\|g_{k+1} \|
\leq \left(\frac{L_H V_{\max}^2}{2\sigma_{k}}\kappa_C^2 \theta^2+\frac{\theta^2 \kappa_B\kappa_C}{\sqrt{\epsilon\sigma_{k}}}+1\right)\|g_k\| .
\eeqn
}
\proof{Let $k \in \calS_k^{curv}$. From \eqref{rhokdef} and \eqref{lowboundcurv},  we obtain that
\[
f(x_k) - f(x_{k+1})
\geq \eta_1 \left(-g_k^\intercal s_k-\frac{1}{2}s_k^\intercal H_ks_k\right)
\geq \frac{\eta_1 \theta^3 \kappa_C^3 }{2 \sqrt{\sigma_{k}}}\|g_k \|^\sfrac{3}{2}. 
\]
Since $\|g_k \| \geq \epsilon$ before termination and that
$\sigma_k \leq \sigma_{\max}$ by Lemma~\ref{boundsigmak}, we obtain
\eqref{threehalfneg}. 

Let us now prove \eqref{succgradnextineq}. 
By using the Lipschitz error bound for the gradient
(\cite[Lemma~2.1]{CartGoulToin20b}), the triangular inequality, the fact that $k \in
\calS_k^{curv}$, \eqref{negcurvvector}, \eqref{negcurvstep},
 and \eqref{skcurvbound},  we obtain that
\begin{align}
\|g_{k+1} \| &\leq \|g_{k+1} - g_k -H_ks_k \| + \|H_ks_k + g_k \| \nonumber \\
&\leq \frac{L_H}{2} \|s_k \|^2  + \|g_k \| +  \|H_ks_k \| \nonumber \\
&= \frac{L_H}{2} \|s_k^{curv} \|^2  + \|g_k \| +  \|H_ks_k^{curv} \| \nonumber\\
&\leq \frac{L_H V_{\max}^2}{2\sigma_k}\kappa_C^2 \theta^2 \|g_k \|  + \|g_k \| +  \|H_ks_k^{curv} \|.\label{gkplus-bound}
\end{align}
Now, using \eqref{negcurvvector}, \eqref{negcurvstep} again,
\begin{align*}
\|H_ks_k^{curv} \|
&= \theta \kappa_C\sqrt{\frac{\|g_k\|}{\sigma_k}}\|H_kV_pu_k\|= \theta \kappa_C\sqrt{\frac{\|g_k\|}{\sigma_k}}\sqrt{u_k^{\intercal}\hH_k^2 u_k}\\
&\leq \theta^2 \kappa_C\sqrt{\frac{\|g_k\|}{\sigma_k}}|\lambda_{\min}(\hH_k)| \leq \theta^2 \kappa_C\sqrt{\frac{\|g_k\|}{\sigma_k}}|\lambda_{\min}(H_k)|.
\end{align*}
Hence \req{gkplus-bound} together with AS.4 and the fact $\|g_k \| \geq
\epsilon$ before termination, give that
\begin{align*}
\|g_{k+1} \|
&\leq \frac{L_H V_{\max}^2}{2 \sigma_k} \kappa_C^2 \theta^2 \|g_k\|+\|g_k \|
+\theta^2 \kappa_B\kappa_C\sqrt{\frac{\|g_k\|}{\sigma_{k}}} \\
&= \left( \frac{L_H V_{\max}^2}{2 \sigma_k} \kappa_C^2 \theta^2
+\frac{ \theta^2 \kappa_B\kappa_C }{\sqrt{\sigma_k\|g_k\|}}+1\right)\|g_k\| \\ 
&\leq \left(  \frac{L_H V_{\max}^2}{2 \sigma_k} \kappa_C^2 \theta^2+\frac{\theta^2 \kappa_B\kappa_C}{\sqrt{\sigma_k\epsilon}}+1\right)\|g_k\|,
\end{align*}
yielding \eqref{succgradnextineq}.
}

\noindent
This lemma is the only result requiring AS.4 or its weaker formulation
\req{boundnegcurv}. Note that this assumption is only required along directions of
negative curvature, which we expect to occur rarely in practice for
suitably large choices of $\kappa_C$. 

After proving a lower bound on the quadratic's decrease when $k \in
\calS_k^{def}$, we now exhibit a relationship between the decrease on
the objective function decrease and gradient both at iteration $k$
and $k+1$ for $k \in \calS_k^{neig} \cup \calS_k^{def}$. This is also
where  the two global conditions \eqref{Vpspacequality} and
\eqref{Vpspacequalitydef} on the subspace $V_p$ will be useful.
Moreover, we also prove an inequality between the norms of the
gradient at two successive iterations, similar to  \eqref{succgradnextineq}. 

\llem{sknlowerbound}{Suppose that AS.1 and AS.3 hold and that
$k \in  \calS_k^{neig} \cup \calS_k^{def}$ before termination.  Then
\beqn{gkplusonenewtonbound}
\|g_{k+1} \|
\leq \left( \frac{L_H V_{\max}^3 (1+\kappa_\theta)}{2 \varsigma_1^2 \sigma_{k}}
      +\frac{2 \kappa_b \sqrt{V_{\max}} }{ \varsigma_1} + \kappa_b \kappa_C \sqrt{V_{\max}} \right)  (1+\kappa_\theta) \|g_k \|
\eeqn		
and
\begin{align}\label{fdecrNewton}
f(x_k) - f(x_{k+1}) &\geq  \eta_1 \, \varsigma_{\min} \sqrt{ \sigkgk}\nonumber \\
&\hspace*{-2em} \left(  \frac{-(2+\kappa_C)  \kappa_b \sqrt{ \kappa_a \sigkgk}
  + \sqrt{ ( \kappa_b (2+\kappa_C))^2{ \kappa_a \sigkgk} + 2\,V_{\max}^2L_H  \|g_{k+1} \|}}  {L_H V_{\max}^2} \right)^2
\end{align}
where 
\beqn{varsigmamindef}
\varsigma_{\min} \eqdef \min \left( \frac{1-2\varsigma_2}{2}, \, 1-\varsigma_3\right).
\eeqn
 }
\proof{
Consider first the case where $k \in \calS_k^{neig}$. By using  the Lipschitz error bound
for the gradient (\cite[Lemma~2.1]{CartGoulToin20b}), that \eqref{Vpspacequality} holds, $r_k^{neig}$ expression \eqref{gkexpressionnewton}, the
condition on $\|r_k^{neig}\|$ \eqref{residcond} and the fact that
$\minlamdavp \leq \kappa_C \sqrt{\sigkgk}$ for $k \in \calS_k^{neig}$,
we deduce that
\begin{align}\label{toreuse}
\|g_{k+1} \| &\leq \|g_{k+1} - H_k s_k^{neig} - g_k \| + \|H_k s_k^{neig} + g_k \| \nonumber \\
&\leq \frac{L_H}{2}\|s_k^{neig}\|^2 + \kappa_b \| \hH_k   y_k^{neig} + \hg_k \| \nonumber \\
&\leq \frac{L_H}{2}\|s_k^{neig}\|^2+ \kappa_b (\sqrt{\sigkgk}+\minlamdavp) \| y_{k}^{neig}\|+ \kappa_b \|r_k^{neig}\| \nonumber \\
&\leq \frac{L_H}{2}\|s_k^{neig}\|^2+ \kappa_b(1+\kappa_C) \sqrt{\sigkgk}\| y_{k}^{neig}\|+ \kappa_b \varsigma_3 \sqrt{\sigkgk }  \| y_{k}^{neig} \|.
\end{align}
Using now \eqref{ypbound} and \eqref{sknnormbound} in the last inequality
\beqn{gkplusonefirstcase}
\|g_{k+1} \|
\leq \left(\frac{L_H V_{\max}^3}{2\sigma_{k}}(1+\kappa_\theta)+\kappa_b(1+\kappa_C) \sqrt{V_{\max}}  + \varsigma_3 \kappa_b \sqrt{V_{\max}} \right) (1+\kappa_\theta)\|g_k\|.
\eeqn
Consider now $k \in \calS_k^{def}$. By arguments similar to those
used for \eqref{toreuse}, this time with \eqref{gkconvexpression}, \eqref{Vpspacequalitydef} and
\eqref{residcondconv}, we obtain that
\begin{align}\label{convcase}
\|g_{k+1} \|
&\leq \|g_{k+1}-H_k s_k^{def}-g_k \| + \|H_k s_k^{def}+g_k \| \nonumber  \\
&\leq \|g_{k+1}-H_k s_k^{def}-g_k \| + \kappa_b \| \hH_k y_{k}^{def} + \hg_k \| \nonumber \\
&\leq \frac{L_H}{2}\|s_k^{def}\|^2+ \kappa_b \sqrt{\kappa_a\sigkgk}\| y_{k}^{def}\|+ \kappa_b \|r_k^{def}\| \nonumber  \\
&\leq \frac{L_H}{2}\|s_k^{def}\|^2+ \kappa_b \sqrt{\kappa_a\sigkgk}\|y_{k}^{def}\|+ \kappa_b \varsigma_2\sqrt{\kappa_a\sigkgk}\|y_k^{def}\|.
\end{align}
Bounding $\|s_k^{def} \|$ with \eqref{skdefbound} and utilizing  \eqref{skconvexbound}  yields that
\beqn{gkplusoneseccase}
\|g_{k+1} \| \
\leq \left( \frac{L_H (1+\kappa_\theta) V_{\max}^2} {2 \varsigma_1^2 \kappa_a \sigma_{k}}
    + \frac{\kappa_b(1+\varsigma_2)}{ \varsigma_1} \right) \left(1+\kappa_\theta \right) \|g_k \|, 
\eeqn
so that taking the larger bound for both \eqref{gkplusonefirstcase}
and \eqref{gkplusoneseccase} and using the bounds $\varsigma_1 < 1$,  $\varsigma_2 <
\frac{1}{2}$, $\varsigma_3 < 1$, $V_{\max} \geq 1$ and $\kappa_b \geq 1$ gives  \eqref{gkplusonenewtonbound}. 

Finally, from \eqref{convcase}, \eqref{toreuse}, \eqref{skneqyp}, \eqref{skdefbound}, the bounds
$\max(\varsigma_3,\varsigma_2) < 1$ and  $\kappa_a \geq 1$, we obtain
that, for $k \in \calS_k^{def} \cup \calS_k^{n}$,
\[
\frac{L_H V_{\max}^2}{2} \|y_{k}\|^2 +  \kappa_b (2+\kappa_C) \sqrt{ \kappa_a \sigkgk} \|y_{k}\| - \|g_{k+1} \| \geq 0.
\]
Hence $\|y_{k}\|$ is larger than the positive root of this quadratic and therefore
\[
\|y_{k}\|
\geq \frac{-\kappa_b (2+\kappa_C) \sqrt{ \kappa_a \sigkgk}
    + \sqrt{ \kappa_b^2(2+\kappa_C)^2{\kappa_a \sigkgk}+2L_H V_{\max}^2 \|g_{k+1}\|}}{L_H V_{\max}^2} > 0.
\]
We then deduce  \eqref{fdecrNewton} from this inequality,
\eqref{rhokdef}, the lower bounds on the quadratic decrease for
$k \in \calS_k^{neig}$ or $k\in\calS_k^{def}$ (\eqref{lowboundnewton} and
\eqref{lowboundconv} respectively) and the definition of
$\varsigma_{\min}$ in \eqref{varsigmamindef}. 
}

\noindent
The bound \eqref{fdecrNewton} is not sufficient for deriving the
required $\mathcal{O}\left(\epsilon^{-3/2}\right)$ optimal complexity
rate because the decrease depends on both $\|g_{k+1} \|$ and $\|g_k
\|$. Indeed, when $\|g_{k+1} \|\ll\|g_k \| $, the right-hand side of
\eqref{fdecrNewton} tends to zero. To circumvent this difficulty, the
next lemma borrows some elements of \cite[Theorem~1]{Misc21} and
partitions $\calS_k^{neig} \cup \calS_k^{def}$ in two further
subsets. The minimum decrease on the objective function is of the
required magnitude in the first one while no meaningful information
can be derived on the decrease on the function value in the second,
albeit the magnitude of the gradient at the next iteration is
halved. The bounds \eqref{gkplusonenewtonbound} and 
\eqref{succgradnextineq} are then used to bound the cardinal of
the latter set.

\llem{halfnextgradstep}{
Suppose that AS.1, AS.3 and AS.4 hold and that 
$\calS_k^{neig} \cup \calS_k^{def} $ is partitioned as
\beqn{divSkn}
\calS_k^{neig} \cup \calS_k^{def}  = \calS_k^{decr} \cup \calS_k^{divgrad}
\eeqn
where 
\beqn{Skndecr}
\calS_k^{decr} \eqdef \{ k \in \calS_k^{neig} \cup \calS_k^{def}, \, \, \sigkgk \leq \kappa_m 2 L_H \|g_{k+1} \|\},
\eeqn 
\beqn{Skndivgrad}
\calS_k^{divgrad} \eqdef \{ k \in \calS_k^{neig} \cup \calS_k^{def} , \, \, \sigkgk > \kappa_m 2 L_H \|g_{k+1} \|\}
\eeqn
with 
\beqn{kappamdef}
\kappa_m \eqdef  {\gamma_3 \max \left( \frac{\sigma_0}{L_H}, \frac{\varsigma_{\max}}{6(1-\eta_2)} \right)}.
\eeqn
Then, for all $k \in \calS_k^{decr}$,
\beqn{lowfuncdecrSknfirst}
f(x_k) - f(x_{k+1})
\geq \frac{\eta_1 \,\varsigma_{\min} (\sigma_{k} \|g_k\| )^\sfrac{3}{2}}
          {\left(\kappa_m  L_H \left( \kappa_b  (2+ \kappa_C) \sqrt{\kappa_a}
       + \sqrt{ (\kappa_b  (2+\kappa_C))^2 \kappa_a + \frac{V_{\max}^2}{\kappa_m}  } \right) \right)^2}.
\eeqn
Moreover,
\beqn{upperSkndivgrad}
|\calS_k^{divgrad} |
\leq \kappa_n |\calS_k^{decr}|+\left(\frac{1}{2\log(2)}|\log(\epsilon)|
     +\kappa_{curv}\right)|\calS_k^{curv}|+\frac{|\log(\epsilon)| + \log(\|g_0\|)}{\log(2)}+1,
\eeqn
where 
\beqn{kappandef}
\kappa_n
\eqdef \frac{1}{\log (2)}\log\left( \frac{L_H(1+\kappa_\theta) V_{\max}^3}{2 \varsigma_1^2\sigma_{\min}} + \frac{2 \sqrt{V_{\max}} \kappa_b}{ \varsigma_1} +
      \sqrt{V_{\max}} \kappa_C \kappa_b\right) + \frac{\log\left( 1+\kappa_\theta\right)}{\log(2)},
\eeqn
\beqn{kappacurvdef}
\kappa_{curv}
\eqdef \frac{1}{ \log(2)} \log \left( \frac{L_H V_{\max}^2}{2 \sigma_{\min}} \kappa_C^2 \theta^2
       + \frac{\theta^2 \kappa_B \kappa_C}{\sqrt{ \sigma_{\min}}} +1 \right).
\eeqn
}
\proof{
Let $k \in \calS_k^{decr}$. Injecting the definition of
$\calS_k^{decr}$ \eqref{Skndecr} in \eqref{fdecrNewton}, we obtain that 
\begin{align*}
f(x_k) - f(x_{k+1})
&\geq \eta_1 \varsigma_{\min} ( \sigkgk)^\sfrac{3}{2} \left(\bigfrac{-(2+\kappa_C) \kappa_b \sqrt{\kappa_a}
 + \sqrt{ (2+\kappa_C)^2 \kappa_b^2 \kappa_a + \frac{V_{\max}^2}{\kappa_m}}}{L_H V_{\max}^2}\right)^2.
\end{align*}
Taking the conjugate both at the denominator and numerator yields \eqref{lowfuncdecrSknfirst}.
	
Let $k \in \calS_k^{divgrad}$. Using the definition of $\kappa_m$ in
\eqref{kappamdef} and that of $\calS_k^{divgrad}$ in
\eqref{Skndivgrad} gives that
\beqn{halfnextgrad}
\|g_{k+1} \|  < \frac{\sigma_k}{\kappa_m L_H} {\frac{\|g_k \|}{2}}
\leq \frac{\sigma_k}{{\gamma_3 \max \left( \frac{\sigma_0}{L_H},\frac{\varsigma_{\max}}{6(1-\eta_2)} \right)}L_H}
       {\frac{\|g_k \|}{2}}
\leq {\frac{\|g_k \|}{2}},
\eeqn
where the last inequality results from the upper bound on $\sigma_{k}$ in \eqref{sigmamax}. 

Successively using the fact that $\calS_k = \calS_k^{decr} \cup
\calS_k^{divgrad} \cup \calS_k^{curv} $, the relationship between
$\|g_{k+1} \|$ and $\|g_k \|$ in the three cases
(\eqref{halfnextgrad}, \eqref{gkplusonenewtonbound} and
\eqref{succgradnextineq}), the fact that $\sigma_{k} \geq
\sigma_{\min}$ in \eqref{gkplusonenewtonbound} and
\eqref{succgradnextineq}, we then deduce that
\begin{align*}
\frac{\epsilon}{\|g_0 \|}\leq\frac{\|g_k \|}{\|g_0 \|} &= \prod_{i \in \calS_k\setminus\{k\}}\frac{\|g_{i+1} \|}{\|g_i \|} \\
& = \prod_{i \in \calS_k^{decr}\setminus\{k\}}  \frac{\|g_{i+1} \|}{\|g_i \|}
    \prod_{i \in \calS_k^{divgrad}\setminus\{k\}} \frac{\|g_{i+1} \|}{\|g_i \|}
    \prod_{i \in \calS_k^{curv}\setminus\{k\}} \frac{\|g_{i+1} \|}{\|g_i \|} \\
&\leq \left[\left(\frac{L_H(1+\kappa_\theta) V_{\max}^3}{2 \varsigma_{1}^2\sigma_{\min}}+ \frac{2 \kappa_b \sqrt{V_{\max}} }{ \varsigma_1} + \kappa_C \kappa_b \sqrt{V_{\max}} \right)
  (1+\kappa_\theta)\right]^{|\calS_k^{decr}\setminus\{k\}|} \times   \\
    & \frac{1}{2^{|\calS_k^{divgrad}\setminus\{k\}|}}\times \hspace*{5mm}\left[ \frac{L_H V_{\max}^2}{2 \sigma_{\min}} \kappa_C^2 \theta^2
      + \frac{\theta^2 \kappa_B \kappa_C}{\sqrt{\epsilon \sigma_{\min}}} +1 \right]^{|\calS_k^{curv} \setminus\{k\}| }.
\end{align*}
Now $\varsigma_1 \leq 1$ and thus both terms in brackets are larger
than one.  Moreover, obviously,
$|\calS_k^{decr}\setminus\{k\}|\le|\calS_k^{decr}|$ and
$|\calS_k^{curv} \setminus\{k\}|\le|\calS_k^{curv}|$, so that
\begin{align*}
\frac{2^{|\calS_k^{divgrad}\setminus\{k\}|}\epsilon}{\|g_0 \|}
&\leq \left[\left(\frac{L_H(1+\kappa_\theta) V_{\max}^3}{2 \varsigma_{1}^2\sigma_{\min}}+ \frac{2 \kappa_b \sqrt{V_{\max} }}{ \varsigma_1} + \sqrt{V_{\max}} \kappa_C \kappa_b \right)
(1+\kappa_\theta)\right]^{|\calS_k^{decr}|} \\
&\times \left[ \frac{L_H V_{\max}^2}{2 \sigma_{\min}} \kappa_C^2 \theta^2
+ \frac{\theta^2 \kappa_B \kappa_C}{\sqrt{\epsilon \sigma_{\min}}} +1
\right]^{|\calS_k^{curv}|}.
\end{align*}
Taking logarithms gives that
\begin{align*}
| \calS_k^{divgrad}\setminus\{k\}| \log(2)
& \leq \log\left[ \left(\frac{L_H(1+\kappa_\theta) V_{\max}^3}{2\varsigma_{1}^2 \sigma_{\min}} + \frac{\kappa_b \sqrt{V_{\max}}}{ \varsigma_{1}}+
       \kappa_C \kappa_b \sqrt{V_{\max}} \right) (1+\kappa_\theta)\right] | \calS_k^{decr}| \\ &+ \log(\|g_0 \|) 
+|\log(\epsilon)|+\log\left[\frac{L_H V_{\max}^2}{2 \sigma_{\min}} \kappa_C^2 \theta^2+
       \frac{\theta^2\kappa_B \kappa_C}{\sqrt{\epsilon \sigma_{\min}}} +1 \right] |\calS_k^{curv} |  . 
\end{align*}
We then obtain \eqref{upperSkndivgrad} with the values of $\kappa_n$
and $\kappa_{curv}$ stated in \eqref{kappandef} and \eqref{kappacurvdef} by dividing this last
inequality by $\log(2)$ and using the facts that
$|\calS_k^{divgrad}\setminus\{k\}|\ge|\calS_k^{divgrad}|-1$
and $\frac{1}{\sqrt{\epsilon}} \geq 1$. 
}

\noindent
Combining the previous lemmas, we are now able to state the complexity
of the \al{AN2C} algorithm. Our theorem statement relies on the
observation that the objective function is evaluated once per
iteration, and its derivatives once per successful iteration.

\lthm{Complexity bound}{Suppose that AS.1- AS.4 hold. Then the
  \tal{AN2C} algorithm requires at most   
\[
| \calS_k | \leq \left( \kappa_\star + \frac{\kap{negdecr}}{2 \log(2)}  | \log (\epsilon)| \right)   \epsilon^{-\sfrac{3}{2}}  + \frac{ |\log(\epsilon)|+ \log(\|g_0 \|) }{\log(2)} + 1
\]
successful iterations and evaluations of the gradient and the Hessian and at most 
\begin{align*}
 \left(1+\frac{|\log\gamma_1|}{\log\gamma_2}\right)& \left[ \left( \kappa_\star +  \frac{\kap{negdecr}}{2 \log(2)} |\log (\epsilon)| \right) \epsilon^{-\sfrac{3}{2}} + \frac{|\log(\epsilon)| + \log(\|g_0\|)}{\log(2)} + 1 \right] \\
&+  \frac{1}{\log \gamma_3} \log \left( \frac{\sigma_{\max}}{\sigma_{0}}\right)
\end{align*}
evaluations of $f$ to produce a vector $x_\epsilon$ such that
$\|g(x_\epsilon) \| \leq \epsilon$, where $\kappa_\star$ is defined by
\beqn{kappastardef}
  \kappa_\star \eqdef \kap{decr}  \left(1+ \kap{n} \right) +  \kap{negdecr} (1+\kappa_{curv}),
\eeqn
with
\beqn{kdecr-def}
\kap{decr}
\eqdef
\frac{\left( L_H \kappa_m ( \sqrt{\kappa_a} \kappa_b  (2+\kappa_C)
	+ \sqrt{ \kappa_a (\kappa_b  (2+\kappa_C))^2 + \frac{V_{\max}^2}{\kappa_m}})\right)^2}
     {\eta_1 \varsigma_{\min} \sigma_{\min}^\sfrac{3}{2}}
\eeqn
and
\beqn{knegdecr-def}
\kap{negdecr}
\eqdef \frac{2 (f(x_0) - f_{\rm low})  \sqrt{\sigma_{\max}}}{\eta_1
  \kappa_C^3 \theta^3},
\eeqn
and where $\kap{n}$ and $\kappa_{curv}$ are given by \eqref{kappandef} and \eqref{kappacurvdef}. 
}

\proof{
First note that we only need to prove an upper bound on $|\calS_k^{decr}| $ and $|
 \calS_k^{curv}|$  to derive a bound on $|\calS_k|$ since
\beqn{calSkcurv}
| \calS_k| = |\calS_k^{decr}  | + | \calS_k^{curv}| + |\calS_k^{divgrad} |
\eeqn
and a bound on $|\calS_k^{divgrad} |$ is given by \eqref{upperSkndivgrad}.
We start by proving an upper bound on $| \calS_k^{curv}|$.  Using AS.2, the lower bound on the decrease of the function values
\eqref{threehalfneg} and that $\sigma_{k}\leq \sigma_{\max}$ as stated in Lemma~\ref{boundsigmak}, we derive that, for $k \in\calS_k^{curv}$, 
\[
f(x_0) - f_{\rm low}
\geq \!\!\!\sum_{i \in \calS_k} f(x_i) - f(x_{i+1}) \geq  \sum_{i \in \calS_k^{curv}} f(x_i) - f(x_{i+1})
\geq |\calS_k^{curv}| \, \frac{\eta_1 \kappa_C^3 \theta^3}{2 \sqrt{\sigma_{\max}}}\,\epsilon^\sfrac{3}{2}
\]
and hence that
\beqn{Skcurvbound}
|\calS_k^{curv}|
\leq \frac{2 (f(x_0) - f_{\rm low})  \sqrt{\sigma_{\max}}}{\eta_1 \kappa_C^3 \theta^3} \,\epsilon^{-\sfrac{3}{2}}
= \kap{negdecr}\,\epsilon^{-\sfrac{3}{2}}.
\eeqn
Similarly for $k \in \calS_k^{decr}$, using AS.2,
\eqref{lowfuncdecrSknfirst}, the fact that $\sigma_{k} \geq
\sigma_{\min}$ and $\|g_{k} \| \geq \epsilon$ before termination
yields that
\[
f(x_0) - f_{\rm low}
\geq \sum_{i \in  \calS_k^{decr}} f(x_i) - f(x_{i+1})
\geq \frac{|\calS_k^{decr}|\eta_1 \,\varsigma_{\min}(\sigma_{\min} \epsilon)^\sfrac{3}{2}}
      {\left(L_H \kappa_m(\sqrt{\kappa_a}(2+\kappa_C) \kappa_b  +\sqrt{\kappa_a ( \kappa_b  (2+\kappa_C))^2+\frac{V_{\max}^2}{\kappa_m}})\right)^2}
\]
where $\kappa_m$ is defined in \eqref{kappamdef}. Rearranging the last
inequality yields that
\beqn{Skndecrbound}
|\calS_k^{decr}|
\leq  \frac{\left( L_H \kappa_m ( \sqrt{\kappa_a} \kappa_b  (2+\kappa_C)
      + \sqrt{ \kappa_a (\kappa_b  (2+\kappa_C))^2 + \frac{V_{\max}^2}{\kappa_m}})\right)^2}
      {\eta_1 \varsigma_{\min} \sigma_{\min}^\sfrac{3}{2}} \, \epsilon^{-\sfrac{3}{2}}
=  \kap{decr} \, \epsilon^{-\sfrac{3}{2}}.
\eeqn
Combining now \eqref{Skcurvbound} and \eqref{Skndecrbound} with the
upper-bound  \eqref{upperSkndivgrad} on  $|\calS_k^{divgrad}|$, we
deduce that
\beqn{Skndivgradbound}
|\calS_k^{divgrad}| \leq \kappa_n \kap{decr} \epsilon^{-\sfrac{3}{2}} + \left( \frac{|\log(\epsilon)|}{2 \log(2)} + \kappa_{curv} \right) \kap{negdecr}   \epsilon^{-\sfrac{3}{2}} + \frac{|\log(\epsilon)| + \log(\| g_0\|)}{\log(2)} + 1.
\eeqn
By summing equations \eqref{Skcurvbound}, \eqref{Skndecrbound}, and \eqref{Skndivgradbound} to bound  $|\calS_k|$ in \eqref{Skcurvbound}, while also isolating the terms based on their different orders with respect to $\epsilon$, we obtain that
\beqn{Skbound}
|\calS_k| \leq \left( \kappa_\star +  \frac{ \kap{negdecr}}{2\log(2)}  |\log (\epsilon)| \right) \epsilon^{-\sfrac{3}{2}} + \frac{|\log(\epsilon)| + \log(\|g_0 \|)}{\log(2)} + 1,
\eeqn
where $\kappa_\star$ is defined in \eqref{kappastardef},
thus proving the first part of the theorem. 
The second part is then deduced  from \eqref{Skbound} combined with Lemma~\ref{SvsU}.
}

\noindent
Regrouping all the problem's dependent constant of the last theorem and keeping the worst
dependency w.r.t $\epsilon$, we derive a  $\calO\left( | \log (\epsilon)|\epsilon^{-3/2} \right)$
complexity order in $\epsilon$ that only differs by the factor
$|\log(\epsilon)|$ from the optimal order for nonconvex second-order methods
\cite{CartGoulToin18a}, a factor which is typically small for practical values
of $\epsilon$. The \al{AN2C} algorithm thus enjoys a better complexity order
than that of past hybrid algorithms \cite{CurtRobi18,LiuLiWangYiYang18,
GoldQuanTrot66} for which the order is $\calO\left(\epsilon^{-2}\right)$. However,
it is marginally worse than that of the more complex second-order linesearch of
\cite{RoyeWrig18} which attains the optimal order. Moreover, we
see in the proof of Theorem~\ref{Complexity bound} that the
$|\log\epsilon|$ term  appears because of \eqref{upperSkndivgrad} and
\eqref{Skcurvbound} and we may hope that the  number of $s_k^{curv}$
iterations is typically much less than its worst-case 
$\mathcal{O}\left(\epsilon^{-3/2}\right)$ in practice. The
trust-region algorithm of  \cite{CurtRobiRoyeWrig21} has the same
total complexity as \al{AN2C} although their method requires only
$\mathcal{O}\left( \epsilon^{-3/2}\right)$ gradient and Hessian calls
whereas our algorithm suffers from an additional $|\log(\epsilon)|$
term.  

\numsection{Finding second-order critical points}\label{so-teaser}

Can the \al{AN2C}  algorithm be strengthened to ensure it will compute
second-order critical points? We show in this section under the same
assumptions as that used for its first-order analysis that approximate second order points can be reached. 

The resulting modified algorithm, which we call \al{SOAN2C} (for 
\al{Second-Order AN2C}) makes extensive use of \al{AN2C}, and is
detailed \vpageref{SOAN2C}.

\algo{SOAN2C}{Second-Order Adaptive Newton with Negative Curvature  (\tal{SOAN2C})}{
 \begin{description}
 \item[Step 0: Initialization ] Identical to $\al{AN2C}[{\rm Step~0}]$ with
   $\epsilon\in(0,1]$ now replaced by $\epsilon_1 \in (0,1]$ and $\epsilon_2 \in (0,1]$.
\item[Step 1: Compute current derivatives]
     Evaluate $g_k$ and $H_k$. Terminate if
     \beqn{SO-termination}
     \|g_k \| \leq \epsilon_1 \tim{and} \lambda_{\min}(H_k) \geq - \epsilon_2.
     \eeqn
\item [Step 2: Compute subspace derivatives ] Form $\hg_k$ and $\hH_k$ as in $\al{AN2C}[{\rm Step~2}].$
\item[Step 3: Step calculation ]
     If $\| g_k \| > \epsilon_1$,
     \beqn{FirstOrderStepOptionnal}
     s_k = s_k^{fo}
     = \regstep(\,\hg_k, \hH_k, V_p, \|g_k \|, \sigma_{k}, \kappa_a, \kappa_b, \varsigma_1, \varsigma_2 \,), \tim{ (Optional). }
     \eeqn
     If $s_k^{fo}$ has been successfully defined, go to Step 4. Else, compute
     \beqn{FirstOrderStep}
     s_k = s_k^{fo} =  \newtstep(\, \hg_k,\hH_k,V_p, \|g_k\|, \sigma_k,\kappa_C,\kappa_b, \kappa_\theta, \varsigma_3, \theta \,).
     \eeqn
      Else ($ \, \|g_k \| \leq \epsilon_1 \,$), compute $u_k$ such that 
     \beqn{negcurvsecondorder}
     g_k^\intercal u_k  \leq 0, \, \,  \|u_k\| = 1
     \tim{and}
     H_k u_k = \lambda_{\min}(H_k) u_k,  \, \, 
     \eeqn
      and set 
     \beqn{negcurvstephess}
     s_k = s_k^{so} \eqdef  \frac{-\lambda_{\min}(H_k)}{\sigma_k} u_k.
     \eeqn
\item[Step 4: Acceptance ratio computation ] Identical to $\al{AN2C}[{\rm Step~5}]$.
\item[Step 5: Regularization parameter update ] Identical to $\al{AN2C}[{\rm Step~6}]$.
\end{description}
}

Prior to reaching an approximate
first-order point, we utilize only the $\regstep$ and $\newtstep$
subroutines to generate tentative steps, hence the '$fo$' (first-order)
superscripts in \eqref{FirstOrderStepOptionnal} and
\eqref{FirstOrderStep}. Similar to Section~2, AS.0 is necessary to
obtain a valid step when $\newtstep$ is invoked. Once an approximate
first-order point is reached, further progress towards second-order
stationarity is obtained by exploiting the negative-curvature
direction \eqref{negcurvsecondorder}-\eqref{negcurvstephess}, thereby
justifying the '$so$' (second-order) superscript.

An upper bound on the  evaluation complexity of the \al{SOAN2C}
algorithm is given by the following theorem.

\lthm{ComplexityboundSecond}{Suppose that AS.1--AS.4 hold. Then the
\al{SOAN2C} algorithm requires at most   
\[
| \calS_k |
\leq \kappa_{\star} \epsilon_1^{-\sfrac{3}{2}} + \kap{so}\epsilon_2^{-3}
+ \frac{|\log(\epsilon_1)|}{2 \log(2)} \kap{negdecr}\epsilon_1^{-\sfrac{3}{2}}
+ \left( \frac{|\log(\epsilon_1)| + \log(\kappa_{gpi})}{\log(2)} + 1 \right) (\kap{so} \epsilon_2^{-3} + 1) 
\]
successful iterations and evaluations of the gradient and the Hessian and at most 
\begin{align*}
  \left(1+\frac{|\log\gamma_1|}{\log\gamma_2}\right)
&\Bigg[\kappa_{\star} \epsilon_1^{-\sfrac{3}{2}} + \kap{so}
  \epsilon_2^{-3} + \frac{|\log(\epsilon_1)|}{2 \log(2)} \kap{negdecr}
  \epsilon_1^{-\sfrac{3}{2}} \\
&+ \left( \frac{|\log(\epsilon_1)| + \log(\kappa_{gpi})}{\log(2)} + 1 \right) (\kap{so} \epsilon_2^{-3} + 1) \Bigg]
+   \frac{1}{\log \gamma_3} \log \left( \frac{\sigma_{\max}}{\sigma_{0}}\right)
\end{align*}
evaluations of $f$ to produce a vector $x_\epsilon$ such that
$\|g(x_\epsilon) \| \leq \epsilon_1$ and $\lambda_{\min}(H_{x_\epsilon}) \geq - \epsilon_2$, where
\beqn{kso-def}
\kap{so} \eqdef \frac{2 \sigma_{\max}^2 (f(x_0) - f_{\rm low})}{\eta_1}
\eeqn
$\kap{gpi}$ is defined in \req{gpibound}
and $\kappa_{\star}$, $\kap{negdecr}$ and $\sigma_{\max}$ (defined by \eqref{kappastardef},
\eqref{knegdecr-def} and \req{sigmamax}, respectively)
depend solely on the problem . 
}

As for Theorem~\ref{Complexity bound}, the bound, in which the
$\epsilon_2^{-3}$ term is likely to dominate, differs from standard
one for second-order algorithms seeking second-order points (in
$\calO(\max(\epsilon_1^{-3/2}, \epsilon_2^{-3})$)
\cite[Theorems~3.3.9 and 3.4.6]{CartGoulToin22} by a (modest) factor
$|\log(\epsilon_1)|$. 

To prove Theorem~\ref{ComplexityboundSecond}, we need to take two main issues into account.
The first is that, because the step may be computed using
\req{FirstOrderStepOptionnal}, \req{FirstOrderStep} but also \req{negcurvstephess},
we need to complete the partition of $|\calS_k|$ by introducing
subsets relevant to this new type of steps.  The second is clearly that
negative curvature information must be exploited in order to guarantee a sufficient
decrease of the objective function when it is discovered close to a
first-order critical point.  This leads to a development which
broadly follows the lines of Section~\ref{complexity-s}, extending the
proofs when necessary to handle the more complicated situation.  The
details of this development are given in appendix.

\numsection{Choosing the subspace}\label{numericstep}

In practice, the algorithm crucially depends on how one chooses the
matrix $V_p$ spanning the iteration-dependent subspace,
and we discuss two options. Each of the choices presented below can
be included in both \al{AN2C} and
\al{SOAN2C}, defined in Section~2 and Section~4, respectively. For
conciseness, we only consider \al{AN2C}. 

\subsection{A full-space variant}
A simple choice of $V_p$ is to consider $V_p \eqdef I_n$, that is the
subspace is in fact the whole space. We note that, in this case,
conditions \eqref{Vpspacequalitydef} or \eqref{Vpspacequality} automatically hold.

We define two variants in this context.  The first is called \al{AN2CER}
(for \al{AN2C Exact} using \regstep) exploits the \regstep\ algorithm in order to limit the need
of possibly costly second-order information. The second, potentially
more costly,  is called \al{AN2CE} and does not use the optional
\regstep\ algorithm, therefore making no attempt to avoid
eigenvalue computations.

These variants may be useful for problems in which
systems \req{convexstep} and \req{steplargegrad} may effectively be
solved (for instance using Cholesky factorizations). As we will see
below, they require on average a single such solution/factorization per
iteration. \al{AN2CER} and \al{AN2CE} may thus be attractive in
the large class of applications for which
off-the-shelf linear solvers are available.
The computation of $\lambda_{\min}(H_k)$ also needs to be
feasible but, due to Algorithm~\regstep, this
occurs only rarely in \al{AN2CER}.

\subsection{A Krylov variant}

When the dimension of the problem grows and factorizations become
impractical, one can turn to exploiting Krylov subspaces, as we now
show. The resulting algorithmic variant
will be called \al{AN2CK}, where \al{K} stands here for Krylov, and is obtained by replacing Steps~3 and 4 of
the \al{AN2C} algorithm by Algorithm~\al{AN2CKStep} \vpageref{AN2CKStep}.
In this variant, the subspace generation and step computation are
combined in order to best exploit the structure of the resulting
subproblem. As is common in Krylov-based methods, we assume the
availability of a 'preconditioner', that is a positive-definite matrix $M_k$
approximating the Hessian $H_k$ in the sense that $M_k^{-1}H_k$ is close
to the identity. For clarity, we ignore the iteration subscript $k$ in what follows.

\algo{AN2CKStep}{\al{AN2CKStep}$(\,g, H, \sigma, M, \kappa_C,  \kappa_b, \theta \,)$ }{
	\begin{description}
		\item[Step 0: Initialization] Set $p=1$, $r_1 = g$,
                  $w_1 = M^{-1}r_1$,
                  $\alpha_1=\sqrt{w_1^\intercal r_1}$ and $z_{0} = 0$.
		\item[Step 1: Form the orthonormal basis] Compute
		\beqn{lanczoquantities1}
                z_p = \frac{r_p}{\alpha_p}, \ms
                v_p = \frac{w_p}{\alpha_p}, \ms
                \delta_p = v_p^\intercal H v_p, \ms
                \eeqn
		\beqn{lanczoquantities2}
                r_{p+1} = H v_p - \delta_p z_p - \alpha_p z_{p-1},\ms
                w_{p+1} = M^{-1} r_{p+1},\ms
		\alpha_{p+1} = \sqrt{w_{p+1}^\intercal r_{p+1}},
		\eeqn
		and define 
		\beqn{Vpdef}
		V_p = (v_1 , v_2, \dots , v_p) \in \Re^{n \times p}.
		\eeqn

		\item[Step 2: Newton step computation] 
		Form the subspace Hessian
		\beqn{Tpcompute}
		T_p \eqdef V_p^\intercal H V_p = \begin{pmatrix}
			\delta_1 & \alpha_2  &    &      &  \\
			\alpha_2       & \delta_2  & \alpha_3 &       &    \\
			&   \ddots      & \ddots         & \ddots    &     \\
			&         &           &     \delta_{p-1}      & \alpha_p   \\
			&                 &  & \alpha_{p}           & \delta_p
		\end{pmatrix}
		\eeqn
		and compute its minimum eigenvalue.  
		
		If $\lambda_{\min}(T_p) \leq - \kappa_C \sqrt{\sigma \|g \|}$, go to Step 4. 
		
		Otherwise, solve 
			\begin{equation}\label{steplargegradlanczos}
			\left( T_p + (\sqrt{\sigma \| g\|} + [-\lambda_{\min}(T_p)]_+ ) I_p \right) y_p
                             = - \alpha_1 e_1.
			\end{equation}
		
		\item[Step 3: Check global quality of the solution] If
		  \beqn{Vpspacequalitylan}
                  \sqrt{\alpha_{p+1}^2 (e_p^{\intercal} y_p)^2 + \|T_p
                    y_p + \alpha_1 e_1\|^2} 
			\leq \kappa_b \|T_p y_p  + \alpha_1 e_1 \|,
			\eeqn
			then return
			\beqn{skneiglanczos}
			s = s^{neig} = V_p y_p.
			\eeqn
		Else increment $p$ by one and go back to Step 1.
		
		\item[Step 4: Eigenvector direction] 
			Compute $u$ such that 
			\beqn{negcurvvectorlanczos}
			e_1^\intercal u  \leq 0, \, \,  \|u \| = 1
			\tim{and}
			u^\intercal T_p u \leq \theta \,\lambda_{\min}(T_p).  \, \, 
			\eeqn
			Return
			\beqn{negcurvsteplanczos}
			s = s^{curv}=  \theta \kappa_C  \sqrt{\frac{\|g \|}{\sigma}}\, V_p u.
			\eeqn
	\end{description}
}

Each iteration of the \al{AN2CKStep} algorithm has a moderate cost
(a few vector assignments, one matrix-vector product and --possibly-- the computation of the
smallest eigenvalue of a tridiagonal matrix, see \cite{COakRokh13} and
the references therein for details).  We observe that
\req{lanczoquantities1}-\req{lanczoquantities2} amounts to using the
standard preconditioned Lanczos process
for building an orthonormal (in the $\langle\cdot,M \cdot \rangle$ inner product) basis $V_p$ of successive Krylov
subspaces generate by the preconditioned gradient and Hessian. We therefore build
on existing theory for this process (see
\cite[Section~5.2]{ConnGoulToin00}, for instance).  We note that the 
use of the full Lanczos basis $V_p$ is only requested at the end of
the process (in \req{skneiglanczos} and \req{negcurvsteplanczos}). As
a consequence two options are available for its detailed
implementation: one can store the Lanczos basis vectors as the
iterations proceed and use them at the end of the step computation, or
one can forget them but re-run the necessary Lanczos process to
re-generate them (as has been done in the \al{GALAHAD} library
\cite{GoulOrbaToin03b} for the \al{GLTR} and \al{GLRT} algorithms for
trust-region and regularization subproblems, respectively). Obviously,
$V_p$ and $T_p$ may be updated incrementally in \req{Vpdef} and
\req{Tpcompute}. When updating $T_p$, it is also easy to check if it
remains positive definite by recurring the pivots of its Cholesky
factorization, which are given by
\[
\pi_1 = \delta_1 \tim{and} \pi_p = \delta_p - \alpha_p^2 / \pi_{p-1} \ms
(p>1).
\]
As long as $\pi_p$ stays positive, it is thus unnecessary to compute
 $\lambda_{\min}(T_p)$ since $[-\lambda_{\min}(T_p)]_+$  is then
identically zero in \req{steplargegradlanczos}. Finally, should a
preconditioner $M$ be unavailable, setting $M=I_n$ is possible, in
which case $w_p$ and $z_p$ can be dispensed of because they are
identical to $r_p$ and $v_p$, respectively.

We now verify that, as stated, Algorithm~\al{AN2CK} is a correct
instantiation of Algorithm~\al{AN2C} (without the optional Step~3).

\lthm{EqLanczosgeneric}{
Suppose that
\beqn{Mconds}
\mu_1 \leq \lambda_{\min}(M)  \tim{and} \lambda_{\max}(M) \leq \mu_2
 \eeqn
for some $\mu_2\geq\mu_1>0$.
Then the definitions and conditions \eqref{negcurvvectorlanczos},
\eqref{Vpspacequalitylan} and \eqref{steplargegradlanczos} of
Algorithm~\al{AN2CKStep} are equivalent to \eqref{negcurvvector},
\eqref{Vpspacequality} (with $\kappa_b$ redefined as
$\max(1,\kappa_b\,\sqrt{\mu_2})$ ) and \eqref{steplargegrad} of
Algorithm~\ref{NewtonEigenStep}, respectively.
Moreover, AS.0 holds and \eqref{Tpcompute} is valid.
}  
\proof{
	If $Z_p$ is the matrix whose columns are $z_1,\ldots, z_p$,
        we deduce from \req{lanczoquantities1} and \req{lanczoquantities2} that
	\beqn{VpHTpV}
	H V_p
        = Z_p T_p + \alpha_{p+1} z_{p+1} e_p^\intercal
        = MV_p T_p + \alpha_{p+1} M v_{p+1} e_p^\intercal.
	\eeqn
	Using that $V_p^{\intercal} M v_{p+1} = 0$ yields
        \eqref{Tpcompute}.
	Note also that as $v_1 = \frac{w_1}{\alpha_1}= M^{-1} z_1$ from
        \eqref{lanczoquantities1} and $V_p^\intercal M V_p = I_p$, 
	\beqn{V1prop}
	V_p^\intercal g = \alpha_1 V_p^\intercal z_1 = \alpha_1 V_p^\intercal M v_1 = \alpha_1 e_1.
	\eeqn
	The last identity with the fact that $T_p = V_p^\intercal H
        V_p$ ensures that \eqref{negcurvvectorlanczos} and
        \eqref{steplargegradlanczos} are reformulations of
        \eqref{negcurvvector} and \eqref{steplargegrad}. 
	We now prove that \eqref{Vpspacequalitylan} implies \req{Vpspacequality}.
	Using \eqref{steplargegradlanczos}, \eqref{V1prop},
        \eqref{VpHTpV}, we obtain that
	\begin{align*}
	H s + g
        &= H V_p y_p + \alpha_1 M v_1  = H V_p y_p + \alpha_1  M V_p e_1 \\
	&= H V_p y_p - M V_p T_p y_p - (\sqrt{\sigma \|g \|} + [-\lambda_{\min}(T_p)]_+ ) M V_p y_p \\
	&= \alpha_{p+1} (e_p^{\intercal} y_p) M v_{p+1} - (\sqrt{\sigma \|g \|} + [-\lambda_{\min}(T_p)]_+ ) M V_p y_p.
	\end{align*}
	Since $V_p^\intercal M V_p = I_p$ and $V_p^\intercal M v_{p+1} = 0$,
        we deduce, using \req{steplargegradlanczos} and \req{Vpspacequalitylan}, that
	\begin{align*}
        \|Hs+g\|^2
        &\leq \lambda_{\max}(M)\,(Hs+g)^\intercal M^{-1}(Hs+g)\\
	&= \lambda_{\max}(M)\big[\alpha_{p+1}^2 (e_p^{\intercal} y_p)^2
           + (\sqrt{\sigma \|g \|} + [-\lambda_{\min}(T_p)]_+ )^2 \|y_p\|^2\big]\\
        &= \lambda_{\max}(M) \big[\alpha_{p+1}^2 (e_p^{\intercal} y_p)^2
          + \|T_p y_p + \alpha_1 e_1\|^2 \big]\\
        &\leq \kappa_b^2\,\lambda_{\max}(M) \|T_p y_p + \alpha_1 e_1\|^2,
	\end{align*}
	and \req{Vpspacequality} follows with the redefined $\kappa_b$.
	We finally verify that AS.0 holds. Because
        \[
        1 = \|M^\half V_p\|\geq  \lambda_{\min}(M^\half)\|V_p\|= \sqrt{\lambda_{\min}(M)}\|V_p\|
        \]
        \eqref{Vpbound} holds with
        $V_{\max}=1/\sqrt{\lambda_{\min}(M)}\leq \mu_1^{-1/2}$, where we
        again used \req{Mconds} to derive the last inequality. Moreover, given that $\kappa_b
        \geq 1$, termination necessarily occurs when $p=n$, $V_n^\intercal M V_n = I_n$,
        $V_n$ spans the whole space and $\alpha_{p+1}=0$ in \req{Vpspacequalitylan}.
}

\noindent
The optional Step~3 of Algorithm~\ref{AN2C} is in fact implicitly
contained in Algorithm~\ref{AN2CKStep} since convexity along the current
step (condition \req{possystem}) is verified at each step of the Lanczos
process by checking the positive-definiteness of $T_p$.

Returning now to the complete sequence of minimization iterates, we
see that, whenever the \al{AN2CK} algorithm is used with
iteration-dependent preconditioners $M_k\neq I_n$, Theorems~\ref{Complexity bound} and \ref{ComplexityboundSecond}
remain valid provided \req{Mconds} holds uniformly for all iterations.

\numsection{Numerical illustration}

%\phil{

We now illustrate the behavior of our proposed algorithms on three
sets of test problems from the freely available 
{\sf OPM} collection\footnote{This collection is a subset of the {\sf
  CUTEest} \cite{GoulOrbaToin15b} collection where the test problems are described in Matlab.} \cite{GratToin21c}. The first set contains 119
small-dimensional problems, the second contains 74 medium-size ones,
while the third contains 59 ``largish'' ones. The list of problems and
their dimensions are listed in Tables~\ref{testprobs-s},
\ref{testprobs-m} and \ref{testprobs-l} in appendix.

\subsection{Using the full-space variants}

We use Matlab implementations of \al{AN2CE} and \al{AN2CER} where the involved linear
systems are solved by using the Matlab sparse Cholesky factorization,
and where we have set
\[
\kappa_C = 10^3,\,
\kappa_a = 50 \tim{(\al{AN2CE})  or} 100 \tim{(\al{AN2CER})},\,
\kappa_\theta = 1,\,
\varsigma_1 = \half,\,
\varsigma_2 = \varsigma_3 = 10^{-10},
\]
\[
\sigma_0 = 1, \,
\sigma_{\min} = 10^{-8},\,
\gamma_1 = \half, \,
\gamma_2 = \gamma_3 = 10, \,
\eta_1  =10^{-4} \tim{and}
\eta_2  =0.95.
\]
The values of $\kappa_C$ and $\kappa_a$ were obtained from a
hyper-parameter search\footnote{Covering the choice $\{10^{30}, 10^8, 10^5, 10^3, 10^2, 10\}$ for
$\kappa_C$ and $\{100, 50, 10\}$ for $\kappa_a$.} on the set of small problems.  The
values of $\varsigma_2$ and $\varsigma_3$ are given here for
consistency, but are irrelevant since factorizations are used to
solve the linear systems. Other parameters values are typical of
regularization algorithms.

We compare \al{AN2CE} and \al{AN2CER} with implementations of
the standard adaptive regularization \al{AR2} and trust-region
\al{TR2M}, two well-regarded methods.  All these algorithms use
quadratic approximations of the objective function (i.e. gradients and
Hessians).  The first three also use the same acceptance thresholds
$\eta_1$ and $\eta_2$ and values of $\gamma_1$, $\gamma_2$ and
$\gamma_3$. The \al{TR2M} methods shrinks the trust-region radius by a
factor $\sqrt{10}$ and expands it by a factor 2 (see
\cite[Section~11.2]{CartGoulToin22} for a discussion of the coherence
of these factors between trust-region and adaptive regularization
methods).  The authors are aware that further method-dependent tuning
would possibly result in improved performance, but the values chosen
here appear to work reasonably well for each method. The step
computation is performed in \al{AR2} following \cite[page~67]{CartGoulToin22}
or \cite{BirgGardMartSantToin17} using an (unpreconditioned)
Lanczos approach while a standard Mor\'{e}-Sorensen
method\footnote{Given that our version of \al{AN2C} uses matrix
factorizations, it seems more natural to compare it with a
Mor\'{e}-Sorensen-based trust-region than to one using truncated
conjugate gradients.} is used in \al{TR2M} (see
\cite[Chapter~9]{CartGoulToin22} for details). For \al{AR2}, the step
computation is terminated as soon as
\beqn{ar2test}
\|g_k+H_ks_k\| \leq \half \theta_{sub} \sigma_k\|s_k\|^2
\eeqn 
which slightly differs from the test
$
\|\nabla_s^1 m_k(s_k)\| \leq \half \theta_{sub} \sigma_k\|s_k\|^2
$
used in \cite[page~65]{CartGoulToin22} and
\cite{BirgGardMartSantToin17} while maintaining the desired
$\calO(\epsilon^{-3/2})$ evaluation complexity bound 
(see \cite{GratToin21} for a justification of \req{ar2test} --including the fact that it
more often allows the pure Newton step to be accepted-- or \cite[page~67]{CartGoulToin22}).
The Moré-Sorensen iterations in \al{TR2M} are terminated as soon
as $\|s_k\| \in [(1-\theta_{sub})\Delta_k,(1+\theta_{sub})\Delta_k]$,
where, in both cases, $\theta_{sub} = 10^{-3}$ for $n\leq 100$ and
$10^{-2}$ for $n> 100$.  All experiments were run on a Dell Precison
computer with Matlab 2022b.

\begin{figure}[htb] 
\centerline{
\includegraphics[width=5.7cm]{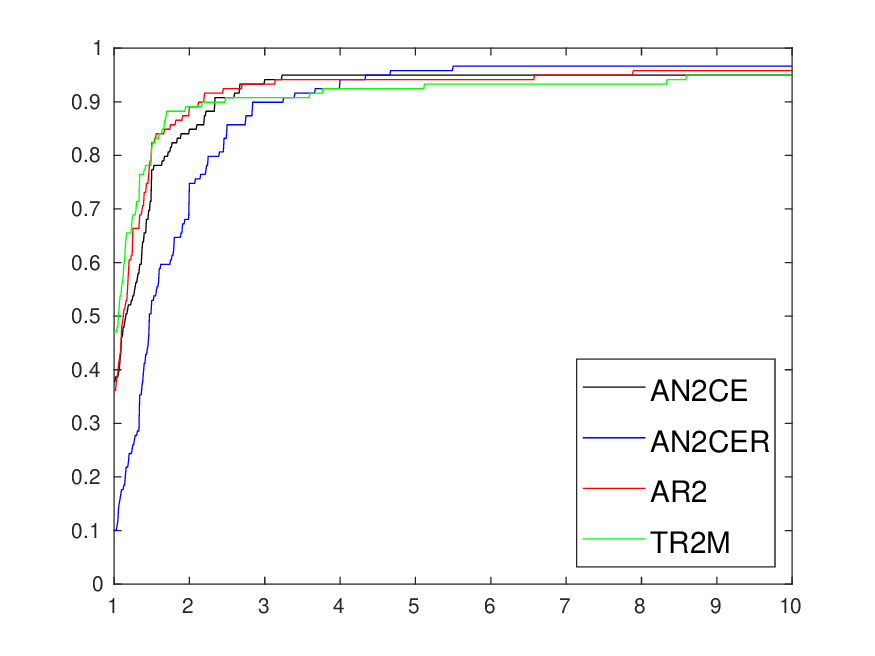}
\hspace*{-5mm}
\includegraphics[width=5.7cm]{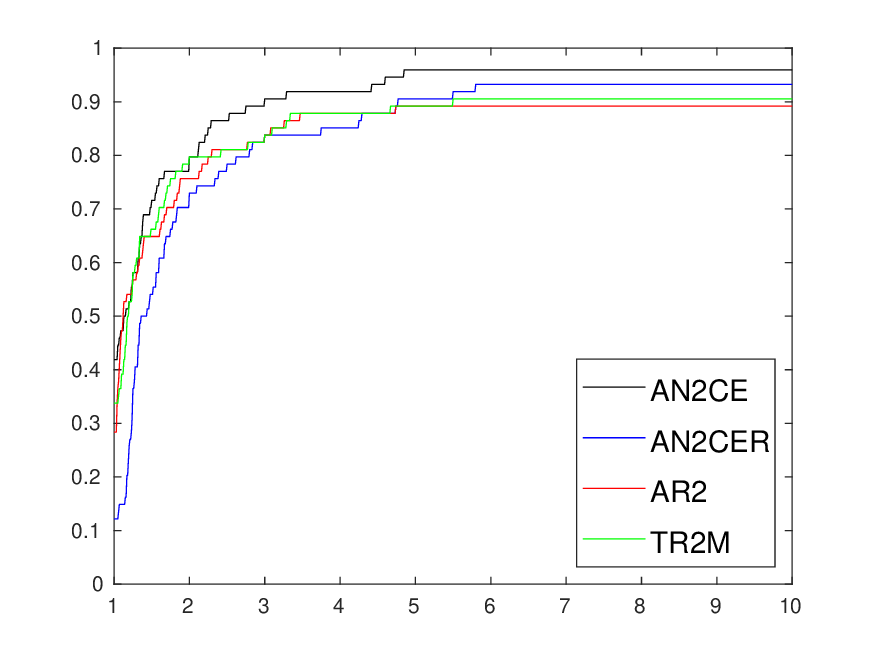}
\hspace*{-5mm}
\includegraphics[width=5.7cm]{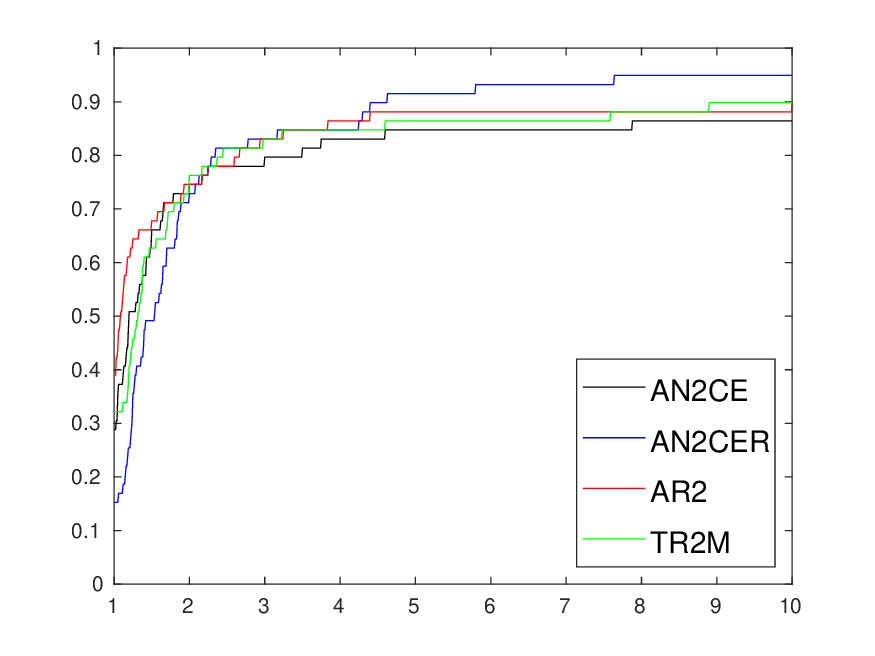}
}
\caption{\label{fig:profile}Full-space variants: iteration performance profiles for 
  {\sf{OPM}} problems (left: small, center: medium, right: largish). We report on the vertical axis
  the proportion of problems for which the number of iterations of each algorithm
  is at most a fraction (given by the horizontal axis) of the smallest
  across all algorithms (see \cite{DolaMoreMuns06}).}
\end{figure}

We discuss our experiments from the efficiency and reliability points
of view. Efficiency is measured, in accordance with the complexity
theory, in number of iterations (or, equivalently, function and
possibly derivatives' evaluations): the fewer the more efficient the 
algorithm. In addition to presenting the now standard performance profiles
\cite{DolaMoreMuns06} for our four algorithms in
Figure~\ref{fig:profile}, we follow \cite{PorcToin17c,GratJeraToin22b}
and consider the derived ``global'' measure $\pi_{\tt algo}$ to be
$\sfrac{1}{10}$ of the area below the curve corresponding to {\tt
  algo} in the performance profile, for abscissas in the interval
$[1,10]$. The larger this area and the closer $\pi_{\tt algo}$ to one,
the closer the curve to the left and top borders of the plot and the
better the global performance.

When reporting reliability, we say that the run of an  
algorithmic variant on a specific test problem is successful if the
gradient norm tolerance $\epsilon = 10^{-6}$ has been achieved in the allotted cpu-time
(1h) and before the maximum number of iterations (5000) is reached.
The $\rho_{\tt algo}$ statistic denotes the percentage of successful
runs taken on all problems in each of the three classes.

\begin{table}[htb]
\begin{center}
\begin{tabular}{|l|r|r|r|r|r|r|}
  \hline
  & \multicolumn{2}{|c|}{small pbs.}& \multicolumn{2}{|c|}{medium pbs.}& \multicolumn{2}{|c|}{largish pbs.}\\
  \hline
{\tt algo} & $\pi_{\tt algo}$  & $\rho_{\tt algo}$ &
$\pi_{\tt algo}$ & $\rho_{\tt algo}$ &$\pi_{\tt algo}$ & $\rho_{\tt algo}$ \\
\hline
\al{AN2CER}  &  0.88 &  96.64  & 0.85  &  93.24 & 0.85 & 94.92 \\
\al{AN2CE}   &  0.91 &  96.64  & 0.91  &  95.95 & 0.81 & 86.44 \\
\al{AR2}     &  0.92 &  97.48  & 0.85  &  93.24 & 0.84 & 93.22 \\
\al{TR2M}    &  0.91 &  94.96  & 0.86  &  93.24 & 0.83 & 91.53 \\
\hline
\end{tabular}
\caption{\label{tab:stats} Efficiency and reliability statistics for
  the {\sf{OPM}} problems (full-space variants)}
\end{center}
\end{table}

Figure~\ref{fig:profile} and Table~\ref{tab:stats} suggest that the
reliability of \al{AN2CE} and \al{AN2CER} is comparable to
that of \al{AR2} and \al{TR2M} for all problem sizes.  They also
indicate that \al{AN2CER} is somewhat slower iteration-wise than
\al{AR2} and \al{TR2M} , but \al{AN2CE} is very comparable. The fact
that the computationally more expensive \al{AN2CE} is often faster than
\al{AN2CER} in terms of iteration numbers is not surprising. Indeed, the
regularization term in \req{steplargegrad} becomes
$\sqrt{\sigma_k \|g_k\|}$ in convex regions, recovering the analysis of
\cite{Misc21,DoikNest23}, whereas \al{AN2CE} regularizes the problem more strongly
in \req{convexstep} (by a factor 10 in our numerical settings) and
therefore may further restricts the steplength. \al{AN2CE} may however
be computationally more intensive\footnote{Most failures of
this algorithm on large problems occurred because the time limit was
reached.} than \al{AN2CER}. Which of the two
algorithms is preferable in practice is likely to depend on the CPU
cost of calculating the Hessian's smallest eigenvalue.

As expected, the call to \newtstep\ in \al{AN2CER} is typically
performed on very few iterations (for less 6.4\% of them for the
small-problems testset) and, when used, results in a
negative-curvature step \req{negcurvstep} even more exceptionally
(less than 1\%). This means in particular
that a single linear-system solve was necessary for approximately 93\%
of all iterations. The \al{AN2CE} variant of course called
\newtstep\ at every iteration, but \req{negcurvstep} was
never actually used.
%Our results also confirm the general effectiveness of a relatively
%simple implementation of the adaptive regularization algorithm
%\al{AR2} using the test \req{ar2test}.

We also ran the \al{SOAN2CE} and \al{SOAN2CER} variants with
$\epsilon_1 = 10^{-6}$ and $\epsilon_2 = 10^{-4}$, but their results are
undistinguishable (for our test sets) from those obtained with
\al{AN2CE} and \al{AN2CER}, except for a final eigenvalue analysis at the
found approximate first-order point, which confirmed in all cases that
the second-order condition \req{SO-termination} did also hold at this
point.  No step of the form \req{negcurvstephess} was ever taken in
our runs, despite the fact that such steps are necessary in theory
(think of starting the minimization at a first-order saddle point).

\subsection{Using the Krylov-based variants}

We ran two variants of the \al{AN2CK} algorithm on our three problem
sets, which differ in how the vector $u$ is chosen in
\req{negcurvvectorlanczos}.  In the first, called \al{AN2CKU}, $u$ is
chosen as the eigenvector associated with the eigenvalue
$\lambda_{\min}(T_p)$.  In the second, called \al{AN2CKYU}, $u$ is
chosen as the sum of the current vector $y_p$ plus a multiple of the
eigenvector associated with $\lambda_{\min}(T_p)$ chosen 
to ensure that the last inequality in \req{negcurvvectorlanczos} holds as
an equality.  An hyper-parameter search on a subset of the
medium-sized test set yielded the values
\[
\kappa_C = 3, \ms \kappa_b = 50 \tim{and} \theta = \half.
\]
None of the tested methods used preconditioning (that is the choice
$M=I_n$ was made throughout). The matrices $V_p$ were stored explicitly.

\begin{figure}[htb] 
\centerline{
\includegraphics[width=5.7cm]{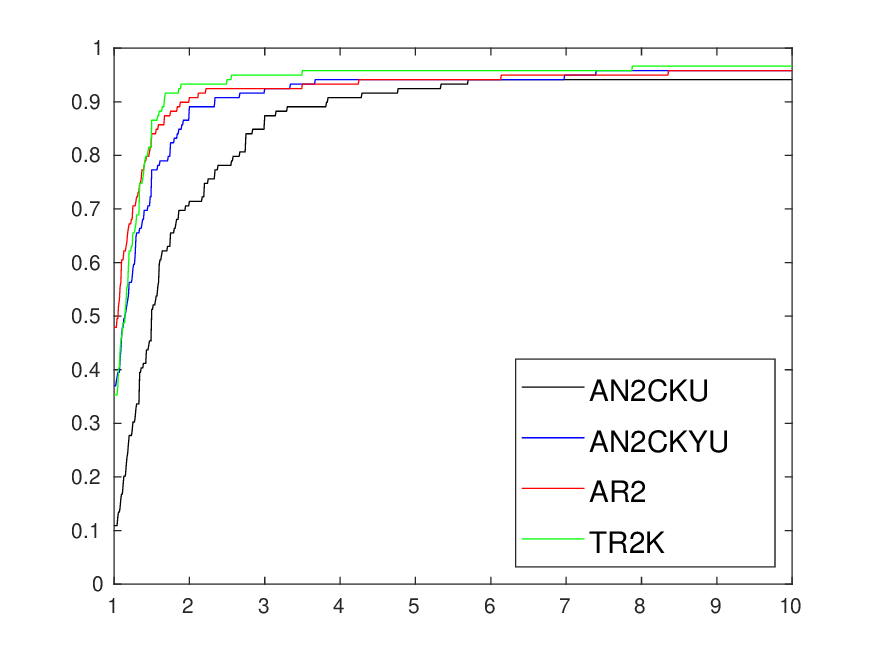}
\hspace*{-5mm}
\includegraphics[width=5.7cm]{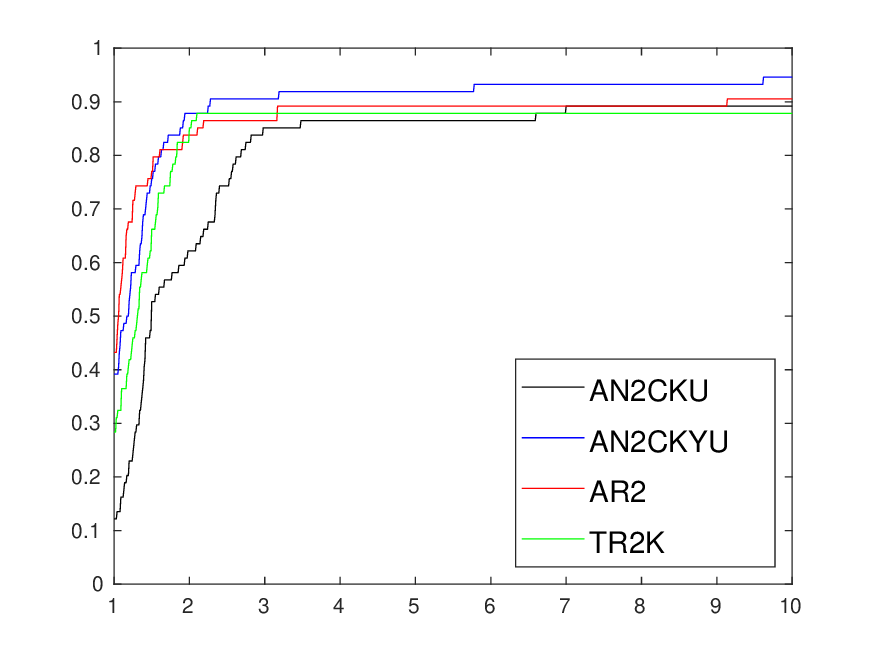}
\hspace*{-5mm}
\includegraphics[width=5.7cm]{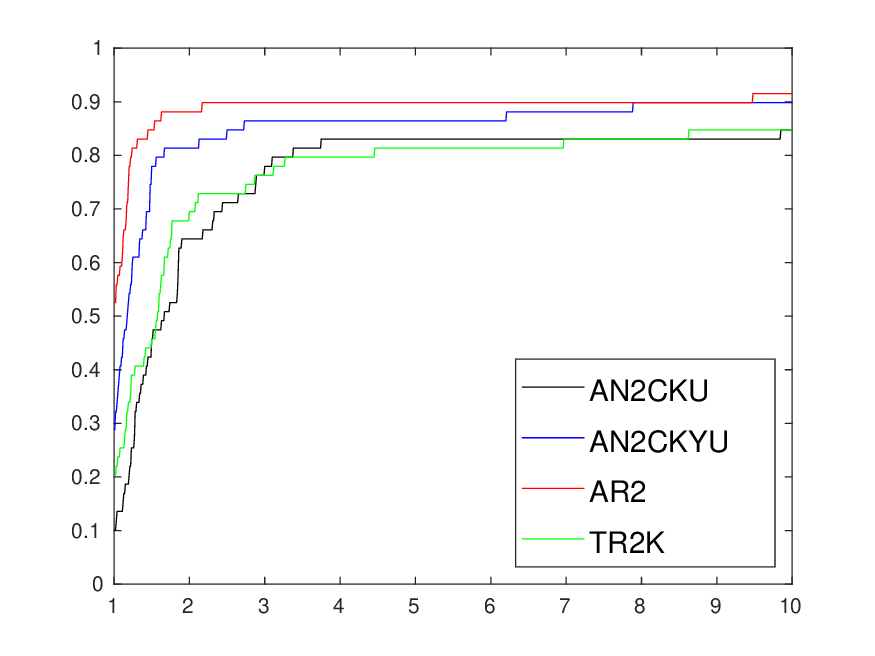}
}
\caption{\label{fig:profile-k}Krylov-space variants: iteration performance profiles for 
  {\sf{OPM}} problems (left: small, center: medium, right: largish). We report on the vertical axis
  the proportion of problems for which the number of iterations of each algorithm
  is at most a fraction (given by the horizontal axis) of the smallest
  across all algorithms (see \cite{DolaMoreMuns06}).}
\end{figure}

We again compared these two variants with \al{AR2} and with \al{TR2K},
an implementation of the trust-region close to \al{TR2M}, but in which the step is computed
by minimizing the quadratic model in the intersection of the trust-region and the successive
Krylov spaces until
\beqn{tr2test}
\|g_k+H_ks_k\| \leq \sfrac{1}{10} \|g_k\|.
\eeqn
The results of our comparison (using the same metrics as in the
previous subsection) are given in Figure~\ref{fig:profile-k} and Table~\ref{tab:stats-k}.

\begin{table}[htb]
\begin{center}
\begin{tabular}{|l|r|r|r|r|r|r|}
  \hline
  & \multicolumn{2}{|c|}{small pbs.}& \multicolumn{2}{|c|}{medium pbs.}& \multicolumn{2}{|c|}{largish pbs.}\\
  \hline
{\tt algo} & $\pi_{\tt algo}$  & $\rho_{\tt algo}$ &
$\pi_{\tt algo}$ & $\rho_{\tt algo}$ &$\pi_{\tt algo}$ & $\rho_{\tt algo}$ \\
\hline
\al{AN2CKU}  &  0.86 &  96.64  & 0.81  &  93.24 & 0.77 & 86.44 \\
\al{AN2CKYU} &  0.91 &  96.64  & 0.90  &  95.95 & 0.85 & 91.53 \\
\al{AR2}     &  0.92 &  97.48  & 0.87  &  93.24 & 0.89 & 93.22 \\
\al{TR2K}    &  0.94 &  96.64  & 0.85  &  87.84 & 0.77 & 84.75 \\
\hline
\end{tabular}
\caption{\label{tab:stats-k} Efficiency and reliability statistics for
  the {\sf{OPM}} problems (Krylov-space variants)}
\end{center}
\end{table}

We observe that \al{AN2CKU} significantly trails the other
variants and is in particular both less efficient and less reliable
than \al{AN2CKYU}, which we explain by the fact that, should a
negative curvature step occur, the former strategy does not exploit
the decrease of the quadratic model already obtained by the ``convex
step'' $y_p$. By contrast, \al{AN2CKYU} appears to be competitive with
both \al{AR2} and \al{TR2K}, irrespective of problem size.

For the \al{AN2CKYU} variant, the average ratio of the number of
matrix-vector products divided by the product of the number of
iterations and the problem size (a ratio which is one if every Lanczos process
takes $n$ iterations) is below 0.5 for small problems, below 0.15 for medium
ones and below 0.03 for large ones.  Negative curvature directions
\req{negcurvsteplanczos} are also used, for this variant, by 0.25\%
of the iterations for small problems, 0.23\% of iterations for medium
ones and never for large ones.

Finally, we also tested \al{SOAN2CKU} and \al{SOAN2CKYU},
the versions of \al{AN2CKU} and \al{AN2CKYU} which enforce second-order
optimality.  As for full-space methods, the results obtained are 
undistinguable from those for \al{AN2CKU} and
\al{AN2CKYU}, except for a final eigenvalue analysis confirming the
approximate second-order optimality of the computed solution.

These early results are encouraging but the authors are aware that
only further experiments will allow a proper assessment of the
method's true potential, both from the number of function/derivatives
evaluations and CPU-usage points of view. Several further algorithmic
developments within the new algorithms are also of interest, including
a possibly better balance between \newtstep\ and \regstep\ in the
full-space version, as well as refinements of the regularization
parameter update \req{sigmakupdate}, possibly in the spirit of
\cite{GoulPorcToin12}.

\numsection{Conclusions and Perspectives}

We have proposed \al{AN2C} and \al{AN2CK}, two second-order
minimization methods for nonconvex problems that alternate, in an
iteration dependent subspace, between Newton and negative-curvature
directions. These methods differ from the more standard trust-region
and adaptive-regularization techniques in that
%, beyond the few requiredeigenvalue calculations,
the involved step computation is free of
further inner iterative processes and only requires the approximate
solution of at most two (but typically one) linear systems per
iteration. We have also proved that these algorithms require at most
$\mathcal{O}\left(|\log(\epsilon)| \epsilon^{-3/2}\right)$ iterations
to obtain an $\epsilon$-approximate first-order critical point.  Our
proof builds on some elements of \cite{Misc21,DoikNest23} for the
convex case and arguments for adaptive regularization
\cite{BirgGardMartSantToin17} and other nonconvex optimization methods
\cite{CurtRobi18,RoyeWrig18}. At each iteration, the algorithms either
take an explicit Newton step or negative curvature when it is
sufficiently large compared to the square root of the gradient. The
norm of the residuals of the Newton step are adjusted dynamically and
different types of solvers can be used to solve the linear systems,
depending on how subspaces are chosen.

An extension of the algorithmic framework ensuring approximate
second-order optimality has also been introduced, and we have proved
that the resulting methods require at most
$\mathcal{O}\left(|\log(\epsilon)|\epsilon^{-3}\right)$ iterations
to achieve its objective.

A first set of numerical experiments with full-space variants
\al{AN2CE} and \al{AN2CER} as well as Krylov-subspaces iterative ones
\al{AN2CKU} and \al{AN2CKYU} indicates that they are very reliable
and competitive with standard techniques in terms of number of
iterations.

The reader may wonder why we haven't considered selecting
iteration-dependent  low-dimensional random subspaces, as has been
advocated in
\cite{CarFowShao22, Shao21} for instance. The main reason is that using the Johnson-Lindenstrauss lemma (the basic tool in such an approach) is possible for defining a probabilistically accurate approximate gradient in the subspace. However, as far as we know, using this idea is problematic for the full Hessian matrix unless it is assumed to be of low rank. We could therefore attempt to follow the Cauchy-point-based
analysis of \cite{CarFowShao22, Shao21}, and hopefully obtain a
probabilistic complexity bound in  $\calO(\epsilon^{-2})$. However, we do not see
at this point how to
design a low-dimensional random-subspace algorithm with an
$\calO(| \log (\epsilon)|\epsilon^{-3/2})$ probabilistic complexity bound for minimizing
functions with general (possibly full-rank) Hessians.

Promising lines for future work include inexact derivatives as in
\cite{YaoXuRoosMaho21,YaoXuRoosMahoWri22}, estimating
the regularization parameter without evaluating the objective function
(as in \cite{GratJeraToin22c}), stochastic variants and the handling
of simple constraints such as bounds on the variables in the spirit of
\cite[Section~14.2]{CartGoulToin22}.

%} %phil

{\footnotesize

  \section*{\footnotesize Acknowlegment}

  The authors are indebted to two excellent referees for their thoughtful
  comments and suggestions.
  
%\bibliographystyle{plain}
%\bibliography{/home/philippe/bibs/refs}

}

\appendix

\appnumsection{Proof of Theorem~\ref{ComplexityboundSecond}}\label{app:so-theory}

As we noted in Section~\ref{so-teaser}, the step in the \al{SOAN2C} algorithm
may be computed using \eqref{FirstOrderStepOptionnal},
\req{FirstOrderStep} or \eqref{negcurvstephess}. The 
notations defining the partition of $|\calS_k|$ remain relevant,
but we complete them by introducing
\[
\calS^{so} \eqdef  \calS \cap \{ s_k = s_k^{so}\},
\ms
\calS_k^{so} \eqdef  \calS_k \cap \{ s_k = s_k^{so}\},
\ms
\calS^{fo} \eqdef  \calS\setminus \calS^{so}
\tim{and}
\calS_k^{fo} \eqdef  \calS_k\setminus \calS_k^{so}.
\]
In addition, for $m\geq \ell \geq 0$, we define
\[
\calS_{\ell,m} \eqdef \calS \cap \iibe{\ell}{m}
\]
and we naturally extend this notation using superscripts identifying
the subsets of $\calS_{\ell,m}$ corresponding to the different  iteration
types identified above.
We  also introduce two index sequences whose purpose is to keep
track of when $s_k= s_k^{fo}$ \eqref{FirstOrderStepOptionnal}-\eqref{FirstOrderStep} or $s_k=s_k^{so}$
\eqref{negcurvstephess} are used, in the sense that
\[
s_k= s_k^{fo} \tim{for} k \in \bigcup_{i\geq0, p_i\geq0} \iibe{p_i}{q_i-1}
\tim{and}
s_k= s_k^{so} \tim{for} k \in \bigcup_{i\geq0} \iibe{q_i}{p_{i+1}-1}.
\]
Formally,
\beqn{q0p0def}
p_0 = \left\{\begin{array}{rl}
  0 &\tim{if } \|g_0 \| > \epsilon_1\\
 -1 &\tim{if } \|g_0 \| \le \epsilon_1,
\end{array}\right.
\tim{ and }
q_0 = \left\{\begin{array}{ll}
  \inf \{  k > 0  \mid \|g_{k} \| \leq \epsilon_1 \} &\tim{if } \|g_0 \| > \epsilon_1\\
 0 &\tim{if } \|g_0 \| \le \epsilon_1.
\end{array}\right.
\eeqn
Then
\beqn{qipidef}
p_i \eqdef \inf \{  k > q_{i-1} \mid \|g_k \| > \epsilon_1\}
\tim{ and }
q_i \eqdef \inf \{  k > p_i \mid \|g_k \| \leq \epsilon_1 \}  \tim{ for } i \geq 1.
\eeqn

\noindent
The following lemma states an important decrease property holding when
\req{negcurvstephess} is used. We also verify that the bound on the
regularization parameter derived in Section~\ref{complexity-s} still applies. 

\llem{sksoproperties}{Suppose that AS.1 and AS.3 hold. Let $k \in \calS^{so}$.
Then 
\beqn{hessdecr}
-g_k^\intercal s_k - \frac {1}{2} s_k^\intercal H_k s_k \geq \frac{1}{2} \sigma_k \|s_k \|^3.
\eeqn
Moreover, the upper bound \eqref{sigmamax} still holds for all $k \geq 0$.
}

\proof{
We obtain from \eqref{negcurvsecondorder} and \eqref{negcurvstephess} that
\[
g_k^\intercal s_k^{so}+\frac{1}{2}(s_k^{so})^\intercal H_k s_k^{so}
\leq \frac{1}{2}  \|s_k^{so}  \|^2 u_k^\intercal H_k u_k 
= \frac{1}{2} \|s_k^{so}  \|^2 \lambda_{\min}(H_k)   
\leq -\frac{1}{2} \sigma_{k} \|s_k^{so}  \|^3,
\]
which gives \eqref{hessdecr}. As in Lemma~\ref{boundsigmak}, we now use
AS.3,  the standard Lipschitz error bound for the
function (see \cite[Lemma~2.1]{CartGoulToin20b}) and \eqref{hessdecr}
to deduce that
\[
1-\rho_k
= \frac{f(x_k+s_k)-f(x_k)-g_k^\intercal s_k-\frac{1}{2}s_k^\intercal H_ks_k}
{-g_k^\intercal s_k - \frac{1}{2} s_k^\intercal H_k s_k}
\leq \frac{L_H\|s_k^{so}\|^3}{6 (\frac{1}{2}\sigma_k\|s_k^{so}\|^3)}
= \frac{L_H}{3 \sigma_{k}}.
\]
Thus, if $\sigma_k \geq \frac{L_H}{3(1-\eta_2)}$, we have that $\rho_k
\geq \eta_2$ and $k$ is a successful iteration.  We may then use the
argument of Lemma~\ref{boundsigmak} and the fact that 
$\varsigma_{\max}$ introduced in \eqref{varsigmamaxdef} is larger than
two as $V_{\max} \geq 1$. Therefore, we deduce that \eqref{sigmamax} also holds for the \al{SOAN2C} algorithm.
}
We now prove an analogue of Lemma~\ref{boundsigmak}, now using the
negative-curvature step as described in
\eqref{negcurvsecondorder}-\eqref{negcurvstephess}. We also bound
the sequence of $\|g_{p_i} \|$. 

\llem{hesscurv}{
Suppose that AS.1, AS.3 and AS.4 hold. Then, for  $k \in \calS^{so}$,
\beqn{decrhess}
f(x_k) - f(x_{k+1}) \
\geq \frac{\eta_1}{2 \sigma_{\max}^2} \epsilon_2^3.
\eeqn
\vspace*{-2mm}
We also have that
\beqn{gpibound}
\|g_{p_i} \| \leq \kappa_{gpi} %g_{\max}
\eqdef\max\left[\|g_0 \|,\left(\frac{L_H \kappa_B^2}{2 \sigma_{\min}^2}+\frac{\kappa_B^2}{\sigma_{\min}}+1\right)\right],
\eeqn 
for all $p_i\geq 0$ as defined in \req{q0p0def}-\req{qipidef}.
}
\proof{
Let $k \in \calS^{so}$. From \eqref{rhokdef} and \eqref{hessdecr},  we obtain that
\[
f(x_k) - f(x_{k+1})
\geq \eta_1 \left(-g_k^\intercal s_k-\frac{1}{2}s_k^\intercal H_ks_k\right)
\geq \frac{\eta_1}{2} \sigma_{k} \|s_k^{so} \|^3. 
\]
Using now that
$\|s_k^{so} \|^3  = \frac{|\lambda_{\min}(H_k)|^3}{\sigma_{k}^3}$
(see \req{negcurvstephess}) in the previous inequality gives that
\[
f(x_k) - f(x_{k+1}) \geq \frac{\eta_1  }{2 {\sigma_{k}}^2} |\lambda_{\min}(H_k)|^3.
\]
Now $|\lambda_{\min}(H_k)| \geq \epsilon_2$  when $s_k^{so}$ is
computed and $\sigma_k \leq \sigma_{\max}$ by Lemma~\ref{sksoproperties},
from which \eqref{decrhess} follows.
Observe now that \req{gpibound} trivially holds if $p_i=p_0=0$.
Consider now $p_i > 0$. From the definition of
$p_i$ and $q_i$ in \eqref{qipidef}, we see that $p_i-1 \in \calS^{so}$. 
Using the Lipschitz error bound for the gradient
(\cite[Lemma~2.1]{CartGoulToin20b}), the triangular inequality,
\eqref{negcurvsecondorder}, \eqref{negcurvstephess},
\eqref{boundnegcurv} (resulting from AS.4),  we obtain that
\begin{align} 
\|g_{p_i} \| &\leq \|g_{p_i} - g_{p_i-1} - H_{p_i-1} s_{p_i-1}^{so} \|  + \|H_{p_i-1} s_{p_i-1}^{so} + g_{p_i-1} \| \nonumber \\
&\leq \frac{L_H}{2} \|s_{p_i-1}^{so} \|^2 + \|g_{p_i-1} \| + \|H_{p_i-1} s_{p_i-1}^{so} \| \nonumber \\
&\leq \frac{L_H | \lambda_{\min}(H_{p_i-1})|^2}{2 \sigma_{p_i-1}^2}
+ \|g_{p_i-1} \| + \frac{|\lambda_{\min}(H_{p_i-1})| \|H_{p_i-1} u_{p_i-1} \| }{\sigma_{p_i-1}} \nonumber \\
&\leq \frac{L_H | \lambda_{\min}(H_{p_i-1})|^2}{2 \sigma_{p_i-1}^2}
       + \|g_{p_i-1} \| + \frac{ |\lambda_{\min}(H_{p_i-1})|^2}{\sigma_{p_i-1}} \nonumber \\
 &= \frac{L_H (- \lambda_{\min}(H_{p_i-1}))^2}{2 \sigma_{p_i-1}^2}
+ \|g_{p_i-1} \| + \frac{ (-\lambda_{\min}(H_{p_i-1}))^2}{\sigma_{p_i-1}} \nonumber \\
&\leq \frac{L_H \kappa_B^2}{2 \sigma_{p_i-1}^2} +\| g_{p_i-1}\| + \frac{\kappa_B^2}{\sigma_{p_i-1}}. \nonumber
\end{align}
But $\|g_{p_i - 1} \| \leq \epsilon_1 \leq 1$  since $p_i - 1 \in
\calS^{so}$ and $\sigma_{k} \geq \sigma_{\min}$ for all $k\geq0$,
which then implies \eqref{gpibound}.
}
In addition to this lemma, all properties of the different steps
derived in Section~3 remain valid because these steps are only
computed for $\|g_k \| > \epsilon_1$. In
particular, \eqref{succgradnextineq} still applies with  $\epsilon =
\epsilon_1$. However, \eqref{upperSkndivgrad} in
Lemma~\ref{halfnextgradstep} may no longer hold because its proof
relies on the fact that $\|g_k \| \geq \epsilon_1$, which is no longer true.
The purpose of the next lemma is to provide an analogue of
\eqref{upperSkndivgrad} for the case where \al{SOAN2C} is used.

\llem{halfnextgradstepsecorder}{
	Suppose that AS.1, AS.3 and AS.4 hold and the \al{SOAN2C}
        algorithm is used. Consider the partition of $\calS_k^{neig} \cup \calS_k^{def}$ into
        $\calS_k^{decr} \cup \calS_k^{divgrad}$ defined in
        Lemma~\ref{halfnextgradstep}  with the same $\kappa_m$
        (defined in \eqref{kappamdef}). 
	Then \req{lowfuncdecrSknfirst} holds for all $k \in \calS_k^{decr}$.
	Moreover,
	\begin{align}\label{Skdivfinal}
		|\calS_k^{divgrad} |&\leq \kappa_n |\calS_k^{decr}|+\left(\frac{1}{2\log(2)}|\log(\epsilon_1)|
	       +\kappa_{curv}\right)|\calS_k^{curv}|
        \nonumber\\
        &\hspace*{5mm}+ \left( \frac{|\log(\epsilon_1)| +
		\log(\kappa_{gpi})}{\log(2)}+ 1 \right) (|\calS_k^{so} | + 1) 
	\end{align}
	where $\kappa_n$ and $\kappa_{curv}$ are defined in
        \req{kappandef} and \eqref{kappacurvdef}  and $\kappa_{gpi}$ is given by \req{gpibound}.
}
\proof{
The proof of \eqref{lowfuncdecrSknfirst} is identical to that used in
Lemma~\ref{halfnextgradstep}.
Moreover, we still obtain \req{halfnextgrad} for $k \in \calS_k^{divgrad}$,
because the definition of $\kappa_m$ in \eqref{kappamdef} is unchanged
and  Lemma~\ref{sksoproperties} ensures that \eqref{sigmamax}
continues to hold for the \al{SOAN2C} algorithm.

We now prove \req{Skdivfinal}.
If  $\calS_k^{fo}$ is empty, then so is its subset $\calS_k^{divgrad}$ and
\req{Skdivfinal} trivially holds.
If $\calS_k^{fo}$ is not empty, we see from the definitions
\eqref{q0p0def}-\eqref{qipidef} that, for some $m\geq0$ depending on $k$,
\beqn{0kdivbiggerepsilon}
\calS_k^{fo}
= \iibe{0}{k} \cap \{ \|g_k \| > \epsilon_1 \}
= \left(\bigcup_{i=0, p_i\geq0}^{m-1} \iibe{p_i}{q_i -1}\right) \cup \iibe{p_{m}}{k}.
\eeqn
Note that the last set in this union is empty unless $k\in\calS^{fo}$,
in which case $p_m \geq 0$. 
Suppose first that the set of indices corresponding to the union in
brackets is non-empty and let $i$ be an index in this set.  Moreover,
suppose also that $p_i < q_i-1$. Using \req{gpibound} and the facts that
$\| g_{q_{i} -1} \| > \epsilon_1$, that the gradient only changes
at successful iterations and that 
$\calS_{p_i,q_i-2}=\calS_{p_i,q_i-2}^{curv}\cup\calS_{p_i,q_i-2}^{divgrad}\cup\calS_{p_i,q_i-2}^{decr}$,
we now derive that 
\begin{align*}
\frac{\epsilon_1}{\kappa_{gpi}} &\leq \frac{\|g_{q_i -1 } \|}{\| g_{p_i}
  \|} = \prod_{j = p_i}^{q_i - 2} \frac{\|g_{j+1} \|}{\|g_j \|} =
\prod_{j \in \calS_{{p_i, q_i - 2}} } \frac{\|g_{j+1} \|}{\|g_j \|}
\\  &= \prod_{j \in \calS_{p_i, q_i - 2}^{decr}}  \frac{\|g_{j+1}
  \|}{\|g_j \|}
\prod_{j \in \calS_{p_i, q_i - 2}^{curv}}  \frac{\|g_{j+1} \|}{\|g_j \|} \prod_{j \in \calS_{p_i, q_i - 2}^{divgrad}}  
\frac{\|g_{j+1} \|}{\|g_j \|}\\
 &\leq \left( \left(\frac{L_H(1+\kappa_\theta) V_{\max}^3}{2 \varsigma_{1}^2
  \sigma_{\min}}+ \frac{2 \kappa_b \sqrt{V_{\max}}}{ \varsigma_1}
+ \kappa_C \kappa_b \sqrt{V_{\max}} \right)(1+\kappa_\theta)\right)^{| \calS_{p_i, q_i-2}^{decr} |}
\times  \\
&\hspace*{5mm} \frac{1}{2^{| \calS_{p_i, q_i - 2}^{divgrad} |}} \times
\left( \frac{L_H V_{\max}^2}{2 \sigma_{\min}} \kappa_C^2 \theta^2
x+ \frac{\theta^2\kappa_B \kappa_C}{\sqrt{\epsilon_1 \sigma_{\min}}} +1 \right)^{|\calS_{p_i, q_i - 2}^{curv} | } 
\end{align*}
where we used  \eqref{gkplusonenewtonbound}, \eqref{succgradnextineq}
and \eqref{halfnextgrad} to derive the last inequality. Rearranging
terms, taking the log, using the inequality 
$|\calS_{p_i,q_i-2}^{divgrad}| \geq |\calS_{p_i,q_i-1}^{divgrad}|-1$ and
dividing by $\log(2)$ then gives that
\[
(|\calS_{p_i,q_i-1}^{divgrad}|-1) + \frac{\log (\epsilon_1)-\log(\kappa_{gpi})}{\log(2)}
\leq \kappa_n |\calS_{p_i, q_i - 2}^{decr}| + \left( \frac{| \log(\epsilon_1)|}{2 \log(2)}
+ \kappa_{curv}\right) |\calS_{p_i, q_i - 2}^{curv}|
\]
with $\kappa_n$ and $\kappa_{curv}$ given by \req{kappandef} and \req{kappacurvdef}.
Further rearranging this inequality and using the fact that
$|\calS_{p_i, q_i - 2}| \leq |\calS_{p_i, q_i - 1}|$ for the different
types of step, we obtain that
\beqn{Skdivpiqi}
|\calS_{p_i, q_i - 1}^{divgrad}|
\leq \kappa_n |\calS_{p_i,q_i-1}^{decr}|
    + \left(\frac{|\log(\epsilon_1)|}{2\log(2)}+\kappa_{curv}\right)|\calS_{p_i,q_i-1}^{curv}|
    + \frac{|\log (\epsilon_1)| + \log(\kappa_{gpi}) }{\log(2)}  + 1.
    \eeqn
If now $p_i = q_i-1$, then clearly $|\calS_{p_i,q_i-1}^{divgrad}|\leq 1$ and \req{Skdivpiqi} also holds.
Using the same reasoning  when $\iibe{p_m}{k}$ is non-empty, we  derive that,
\beqn{Skdivpmk}
|\calS_{p_m, k}^{divgrad}|
\leq \kappa_n |\calS_{p_m, k}^{decr}|
     + \left(\frac{|\log(\epsilon_1)|}{2\log(2)}+\kappa_{curv}\right)|\calS_{p_m,k}^{curv}|
     +  \frac{|\log (\epsilon_1)| + \log(\kappa_{gpi})}{\log(2)}  + 1,
 \eeqn
and this inequality also holds if  $\iibe{p_m}{k}=\emptyset$ since $\calS_{p_m, k}^{divgrad}\subseteq\iibe{p_m}{k}$.
Adding now \eqref{Skdivpiqi} for $i \in \iiz{m}$ and
\req{Skdivpmk} to take \req{0kdivbiggerepsilon} into account gives that
\[%beqn{Skdivk}
|\calS_{ k}^{divgrad}|
\leq \kappa_n |\calS_{ k}^{decr}|
     + \left(\frac{| \log (\epsilon_1)|}{2\log(2)}+\kappa_{curv}\right) |\calS_{k}^{curv}|
     + \left(\frac{|\log (\epsilon_1)| + \log(\kappa_{gpi}) }{\log(2)}+1\right) (m+1).
\]%eeqn
As \eqref{0kdivbiggerepsilon} divides $\calS_{k}^{fo}$
  into $m+1$ consecutive sequences, these sequences are then separated
  by at least a second-order step, so that $m \leq \calS_{ k}^{so}$
  and \req{Skdivfinal} follows. 
}

\noindent
Equipped with this last lemma and the results of
Sections~\ref{thealgo-s} and \ref{complexity-s}, we
may finally establish the worst-case iteration/evaluation complexity
of the \al{SOAN2C} algorithm and prove Theorem~\ref{ComplexityboundSecond} itself.
\vspace*{3mm}
\proof{Note that the bounds \eqref{Skcurvbound} and
\eqref{Skndecrbound} derived in the proof of
Theorem~\ref{Complexity bound} are still valid because they only cover
steps computed using \al{AN2C}, so that we now need to focus on bounding $\calS_k^{so}$.
 Using AS.2
and the lower bound on the decrease of the function values
\eqref{decrhess}, we derive that, for $k \in\calS^{so}$, 
\[
f(x_0) - f_{\rm low}
\geq \!\!\sum_{i \in \calS_k} f(x_i) - f(x_{i+1}) \geq  \sum_{i \in \calS_k^{so}} f(x_i) - f(x_{i+1})
\geq |\calS_k^{so}| \, \frac{\eta_1}{2 \sigma_{\max}^2}\, \epsilon_2^3,
\]
and therefore that
\beqn{Skhessbound}
|\calS_k^{so}| \leq  \frac{2 \sigma_{\max}^2 (f(x_0) - f_{\rm low})}{\eta_1} \, \epsilon_2^{-3} = \kap{so} \epsilon_2^{-3}.
\eeqn
Injecting now \eqref{Skhessbound}, \eqref{Skndecrbound} and
\eqref{Skcurvbound} in the bound \eqref{Skdivfinal} on $\calS_k^{divgrad}$ yields that
\begin{align*}
|\calS_{ k}^{divgrad}| &\leq \kappa_n \kap{decr} \epsilon_1^\sfrac{-3}{2} + \left( \frac{|\log(\epsilon_1)|}{2 \log(2)} + \kappa_{curv} \right) \kap{negdecr} \epsilon_1^\sfrac{-3}{2}\\  &\hspace*{1em}+ \left( \frac{|\log(\epsilon_1)| + \log(\kappa_{gpi})}{\log(2)} + 1 \right) (\kap{so} \epsilon_2^{-3} + 1). 
\end{align*}

Combining the last inequality with  \eqref{Skhessbound},
\eqref{Skndecrbound} and \eqref{Skcurvbound} in
$|\calS_k| = |\calS_k^{divgrad}| + |\calS_k^{curv}| + |\calS_k^{so}| + |\calS_k^{decr}|$
and the definition of \eqref{kappastardef} gives that
\[
|\calS_{ k} | \leq \kappa_{\star} \epsilon_1^{\sfrac{-3}{2}} + \kap{so} \epsilon_2^{-3} + \frac{|\log(\epsilon_1)|}{2 \log(2)} \kap{negdecr} \epsilon_1^\sfrac{-3}{2} + \left( \frac{|\log(\epsilon_1)| + \log(\kappa_{gpi})}{\log(2)} + 1 \right) (\kap{so} \epsilon_2^{-3} + 1). 
\]
This proves the first part of the theorem. 
The second part follows from the last inequality and Lemma~\ref{SvsU}.
}

\noindent
The factor $|\log(\epsilon_1)|$ by which the bound of
Theorem~\ref{ComplexityboundSecond} differs from
$\calO(\max(\epsilon_1^{-3/2}, \epsilon_2^{-3})$ occurs as a
consequence of \eqref{Skdivfinal}, \eqref{Skhessbound} and
\eqref{Skcurvbound} and one expects that, in practice,
$\eqref{Skhessbound}$ is smaller than
$\mathcal{O}\left( \epsilon_2^{-3}\right)$ so that Newton steps are
taken most often.

\newpage
\appnumsection{The test problems and their dimensions}\label{app:details}

\begin{table}[htb]\tiny
\begin{center}
\begin{tabular}{|l|r|l|r|l|r|l|r|l|r|l|r|}
\hline
Problem & $n$ & Problem & $n$ & Problem & $n$ & Problem & $n$ & Problem & $n$ & Problem & $n$ \\
\hline
argauss       &  3 & chebyqad    & 10 & dixmaanl    & 12 & heart8ls   &  8 & msqrtals    & 16 & scurly10    & 10 \\
arglina       & 10 & cliff       &  2 & dixon       & 10 & helix      &  3 & msqrtbls    & 16 & scosine     & 10 \\ 
arglinb       & 10 & clplatea    & 16 & dqartic     & 10 & hilbert    & 10 & morebv      & 12 & sisser      &  2 \\
arglinc       & 10 & clplateb    & 16 & edensch     & 10 & himln3     &  2 & nlminsurf   & 16 & spmsqrt     & 10 \\
argtrig       & 10 & clustr      &  2 & eg2         & 10 & himm25     &  2 & nondquar    & 10 & tcontact    & 49 \\
arwhead       & 10 & cosine      & 10 & eg2s        & 10 & himm27     &  2 & nzf1        & 13 & tquartic    & 10 \\
bard          &  3 & crglvy      &  4 & eigenals    & 12 & himm28     &  2 & osbornea    &  5 & trigger     &  7 \\
bdarwhd       & 10 & cube        &  2 & eigenbls    & 12 & himm29     &  2 & osborneb    & 11 & tridia      & 10 \\
beale         &  2 & curly10     & 10 & eigencls    & 12 & himm30     &  3 & penalty1    & 10 & tlminsurfx  & 16 \\
biggs5        &  5 & dixmaana    & 12 & engval1     & 10 & himm32     &  4 & penalty2    & 10 & tnlminsurfx & 16 \\
biggs6        &  6 & dixmaanb    & 12 & engval2     &  3 & himm33     &  2 & penalty3    & 10 & vardim      & 10 \\
brownden      &  4 & dixmaanc    & 12 & expfit      &  2 & hypcir     &  2 & powellbs    &  2 & vibrbeam    &  8 \\
booth         &  2 & dixmaand    & 12 & extrosnb    & 10 & indef      & 10 & powellsg    & 12 & watson      & 12 \\
box3          &  3 & dixmaane    & 12 & fminsurf    & 16 & integreq   & 10 & powellsq    &  2 & wmsqrtals   & 16 \\
brkmcc        &  2 & dixmaanf    & 12 & freuroth    &  4 & jensmp     &  2 & powr        & 10 & wmsqrtbls   & 16 \\
brownal       & 10 & dixmaang    & 12 & genhumps    &  5 & kowosb     &  4 & recipe      &  2 & woods       & 12 \\
brownbs       &  2 & dixmaanh    & 12 & gottfr      &  2 & lminsurf   & 16 & rosenbr     & 10 & yfitu       &  3 \\
broyden3d     & 10 & dixmaani    & 12 & gulf        &  4 & mancino    & 10 & s308        &  2 & zangwill2   &  2 \\
broydenbd     & 10 & dixmaanj    & 12 & hairy       &  2 & mexhat     &  2 & sensors     & 10 & zangwill3   &  3 \\
chandheu      & 10 & dixmaank    & 12 & heart6ls    &  6 & meyer3     &  3 & schmvett    &  3 &             &    \\
\hline
\end{tabular}
\caption{\label{testprobs-s} The {\sf OPM} small test problems and their dimension}
\end{center}
\end{table}

\begin{table}[htb]\tiny
\begin{center}
 \hspace*{-5mm}\ \begin{tabular}{|l|r|l|r|l|r|l|r|l|r|l|r|}
    \hline
Problem & $n$ & Problem & $n$ & Problem & $n$ & Problem & $n$ & Problem & $n$ & Problem & $n$ \\
\hline
arglina      & 400 & crglvy    & 400 & dixmaanj  & 600 & fminsurf & 400 &  ncb20c    & 500 & tcontact    & 400 \\  
arglinb      &  50 & cube      & 500 & dixmaank  & 600 & freuroth & 500 &  nlminsurf & 400 & tquartic    & 500 \\  
arglinc      &  50 & curly10   & 500 & dixmaanl  & 600 & helix    & 500 &  nondquar  & 500 & tridia      & 500 \\  
argtrig      &  50 & deconvu   &  51 & dixon     & 500 & hilbert  & 500 &  nzf1      & 520 & tlminsurfx  & 400 \\  
arwhead      & 500 & dixmaana  & 600 & dqrtic    & 500 & hydc20ls &  99 &  penalty1  & 500 & tnlminsurfx & 400 \\  
bdarwhd      & 500 & dixmaanb  & 600 & edensch   & 500 & indef    & 500 &  penalty2  & 100 & vardim      & 500 \\  
brownal      & 500 & dixmaanc  & 600 & eg2       & 400 & integreq & 500 &  penalty3  & 500 & wmsqrtals   & 400 \\  
broyden3d    & 500 & dixmaand  & 600 & eg2s      & 400 & lminsurf & 400 &  powellsg  & 500 & wmsqrtbls   & 400 \\  
broydenbd    & 500 & dixmaane  & 600 & eigenals  & 110 & msqrtals & 400 &  powr      & 500 & woods       & 500 \\  
chandheu     & 500 & dixmaanf  & 600 & eigenbls  & 110 & msqrtbls & 400 &  rosenbr   & 100 &             &     \\  
chebyqad     & 150 & dixmaang  & 600 & eigencls  & 110 & morebv   & 500 &  sensors   & 100 &             &     \\  
clplatea     & 400 & dixmaanh  & 600 & engval1   & 500 & ncb20    & 500 &  scosine   & 500 &             &     \\  
clplateb     & 400 & dixmaani  & 600 & extrosnb  & 500 & ncb20b   & 500 &  spmsqrt   & 997 &             &     \\  
\hline

\end{tabular}
\caption{\label{testprobs-m} The {\sf OPM} medium-size test problems and their dimension}
\end{center}
\end{table}

\begin{table}[htb]\tiny
\begin{center}
 \hspace*{-5mm}\ \begin{tabular}{|l|r|l|r|l|r|l|r|l|r|l|r|}
    \hline
Problem & $n$ & Problem & $n$ & Problem & $n$ & Problem & $n$ & Problem & $n$ \\
\hline
arwhead   & 2000 & dixmaand & 2400 & eg2      & 1600 & integreq  & 2000 & powellsg    & 2000 \\
bdarwhd   & 2000 & dixmaane & 2400 & eg2s     & 1600 & lminsurf  & 4900 & powr        & 2000 \\
broyden3d & 2000 & dixmaanf & 2400 & eigenals & 2550 & msqrtals  & 1600 & rosenbr     & 2000 \\
broydenbd & 2000 & dixmaang & 2400 & eigenbls & 2550 & msqrtbls  & 1600 & spmsqrt     & 1498 \\
clplatea  & 4900 & dixmaanh & 2400 & eigencls & 2550 & morebv    & 5000 & tcontact    & 4900 \\
clplateb  & 4800 & dixmaani & 2400 & engval1  & 2000 & ncb20b    & 2000 & tquartic    & 2000 \\
crglvy    & 4000 & dixmaanj & 2400 & extrosnb & 2000 & ncb20c    & 2000 & tridia      & 2000 \\
cube      & 2000 & dixmaank & 2400 & fminsurf & 4900 & nlminsurf & 4900 & tlminsurfx  & 4900 \\
curly10   & 1000 & dixmaanl & 2400 & freuroth & 2000 & nondquar  & 2000 & tnlminsurfx & 4900 \\
dixmaana  & 2400 & dixon    & 2000 & helix    & 2000 & nzf1      & 2600 & vardim      & 2000 \\
dixmaanb  & 2400 & dqrtic   & 2000 & hilbert  & 2000 & penalty1  & 2000 & woods       & 2000 \\
dixmaanc  & 2400 & edensch  & 2000 & indef    & 2000 & penalty3  & 2000 &             &      \\
\hline
\end{tabular}
\caption{\label{testprobs-l} The {\sf OPM} largish test problems and their dimension}
\end{center}
\end{table}

\newpage

\end{document}